\documentclass[11pt]{article}
\usepackage{amsmath,amssymb,fullpage}

\usepackage{enumerate}

\usepackage{graphicx}

\usepackage{makeidx}
\usepackage[utf8]{inputenc}
\usepackage{wrapfig}
\usepackage{amsmath}
\usepackage{amsfonts}
\usepackage{mathrsfs}
\usepackage{amssymb}
\usepackage{latexsym}
\usepackage{amsbsy}
\usepackage{color}
\usepackage{bbm}

\usepackage[colorlinks,citecolor=blue,urlcolor=blue]{hyperref}

\usepackage{mathrsfs}

\usepackage{wasysym}

\providecommand{\otherindexspace}[1]{}

\usepackage{amsthm}

\newtheorem{theorem}{Theorem}[section]

\newtheorem{lemma}[theorem]{Lemma}
\newtheorem{proposition}[theorem]{Proposition}
\newtheorem{remark}[theorem]{Remark}

\newtheorem{definition}[theorem]{Definition}

\newtheorem{corollary}[theorem]{Corollary}

\newtheorem{assumption}[theorem]{Assumption}

\numberwithin{equation}{section}

\newcommand{\dproof}{\noindent {Proof.} \quad}
\newcommand{\fproof}{\hfill $\square$ \bigskip}

\def\R{{\bf R}}

\def\p{\partial}

\def\cal#1{\mathcal{#1}}

\def \H{\mathbb {H}}
\def \R{\mathbb {R}}

\def\dd{\displaystyle}

\begin{document}

\title{Mixed generalized Dynkin game and stochastic control in a Markovian framework}

\author{Roxana Dumitrescu
\thanks{
Institut für Mathematik, Humboldt-Universität zu Berlin, Unter den Linden 6, 10099 Berlin, Germany,
 email: \textbf{roxana@ceremade.dauphine.fr}. 
 The research leading to these results has received funding from the Region Ile-de-France.}
 \and Marie-Claire Quenez 
 \thanks{LPMA,Universit\'e Paris 7 Denis Diderot, Boite courrier 7012, 75251 Paris cedex 05, France, email: \textbf{quenez@math.univ-paris-diderot.fr}} 
\and  Agn\`es Sulem
\thanks{ INRIA Paris, 3 rue Simone Iff, CS 42112, 75589 Paris Cedex 12, France, and Universit\'e Paris-Est, email: \textbf{agnes.sulem@inria.fr}}}

\date{\today}

\maketitle

\begin{abstract}
We introduce a  mixed {\em generalized} Dynkin game/stochastic control with ${\cal E}^f$-expec\-tation in a Markovian framework. 
We study both  the case when the terminal reward function is  Borelian only and when it is continuous. By using the characterization of the value function of a {\em generalized} Dynkin game via an associated doubly reflected BSDEs (DRBSDE) first provided 
in \cite{DQS2}, we obtain that the value function of our  problem coincides with the value function of an optimization problem for DRBSDEs.
Using this property, we establish a weak dynamic programming principle  by extending some  results recently provided in \cite{DQS3}.
We then show a strong dynamic programming principle in the continuous case, 
which cannot be derived from the weak one. In particular, we have to prove that the value function of the problem is continuous with respect to time $t$, which requires some technical tools of stochastic analysis and  new results on DRBSDEs. 
We finally study the links between our mixed problem and generalized Hamilton--Jacobi--Bellman variational inequalities in both cases.

 \end{abstract}



\section{Introduction}

In this paper, we study a new mixed stochastic control/optimal stopping game in a Markovian framework which can be formulated as follows. We consider two actors $A$ and $B$. 
 Actor $A$, called ``controller/stopper", can control a state process $X^\alpha$ through the selection of a control process $\alpha$, which impacts both the drift and the volatility, and can also choose the duration of the ``game" via  a stopping time  $\tau$.
 Actor $B$, called ``stopper", can only decide when to stop the game via an another stopping time $\sigma$. 
 
 We denote by $\mathcal{T}$ the set of stopping times with values in $[0,T]$, where $T>0$ is a fixed terminal time, 
 and by $\mathcal{A}$ the set of admissible control processes $\alpha$.  For each $\alpha \in \mathcal{A}$,  let
$(X_s^\alpha)$ be a jump diffusion process 
of the form 
$$X_s^\alpha=x+\int_0^sb(X_u^\alpha,\alpha_u)du+\int_0^s\sigma(X_u^\alpha,\alpha_u)dW_u+\int_0^s \int_{{\bf E}}\beta(X_u^\alpha,\alpha_u,e)\Tilde{N}(du,de).$$

 If  $A$  chooses a strategy $(\alpha, \tau)$ and $B$  chooses a stopping time $\sigma$, the associated cost (or gain), is  defined by 
 $g(X_{T}^{\alpha})$ if $A$ and $B$ decide to stop at terminal time $T$,
 $h_1( X_{\tau}^{\alpha})$ if $A$ stops before $B$, and $h_2( X_{\sigma}^{\alpha})$ otherwise.
More precisely, the cost is given by
 \begin{equation}\label{Perf}
 I^\alpha(\tau,\sigma):= 
  h_1( X_{\tau}^{\alpha})\textbf{1}_{ \tau \leq \sigma, \tau <T }+
h_2( X_{\sigma}^{\alpha})\textbf{1}_{\sigma < \tau }+
g(X_{T}^{\alpha})\textbf{1}_{\tau = \sigma=T}.
\end{equation}

 If the strategy of $A$ is given by $(\alpha, \tau)$, 
 the stopper $B$ wants to choose $\sigma$ in order to minimize the   expected cost evaluated under a nonlinear expectation. This nonlinear expectation denoted by $\mathcal{E}^{{\alpha}}$ is defined via a BSDE with jumps with a driver $f^{\alpha}$ which
may depend on the control $\alpha$. 
The minimal   expected cost  for $B$ is then given by  $\inf_{\sigma \in \mathcal{T}} \mathcal{E}_{0, \tau \wedge \sigma}^{\alpha} [ I^\alpha(\tau,\sigma)]$. 
 The aim of  actor $A$ is to maximize this quantity over all  choices of $(\alpha,\tau)$, which leads to 
the following mixed optimization problem:
\begin{equation}\label{Value}
\sup_{(\alpha, \tau) \in \mathcal{A} \times \mathcal{T}} 
\inf_{\sigma \in \mathcal{T}} \mathcal{E}_{0, \tau \wedge \sigma}^{\alpha} [ I^\alpha(\tau,\sigma)].
\end{equation}

In the special case when the first player can only act on the duration $\tau$ of the game (i.e. when there is no control $\alpha$), this problem reduces to a {\em generalized} Dynkin game 
that we have introduced in \cite{DQS2}. It is proved there that the value function at $0$, given by 
\begin{equation}\label{optionjeu}
\sup_{\tau \in \mathcal{T} }
\inf_{\sigma \in \mathcal{T}}\mathcal{E}_{0, \tau \wedge \sigma} [ h_1( X_{\tau})\textbf{1}_{ \tau \leq \sigma, \tau <T }+
h_2( X_{\sigma})\textbf{1}_{\sigma < \tau }+
g(X_{T})\textbf{1}_{\tau = \sigma=T}],
\end{equation}

is characterized via a doubly reflected BSDE. 

Note that Problem \eqref{Value} can be seen as a mixed {\em generalized} Dynkin game/Stochastic control problem since 
\begin{equation}\label{Value2}
\sup_{(\alpha, \tau) \in \mathcal{A} \times \mathcal{T}} 
\inf_{\sigma \in \mathcal{T}} \mathcal{E}_{0, \tau \wedge \sigma}^{\alpha} [ I^\alpha(\tau,\sigma)]
= \sup_{\alpha \in \mathcal{A} } (\,\sup_{\tau \in \mathcal{T}}
\inf_{\sigma \in \mathcal{T}}\mathcal{E}_{0, \tau \wedge \sigma}^{\alpha} [ I^\alpha(\tau,\sigma)]\,).
\end{equation}
The control $\alpha$ can then be interpreted as an ambiguity parameter on the model which affects both the drift and the volatility of the underlying state process. 

Using the characterization of the solution of a DRBSDE as the value function of a {\em generalized} Dynkin game we have provided in  \cite{DQS2}, Problem \eqref{Value2} corresponds to  an optimization problem on DRBSDEs with jumps.

\

In this paper, we first give
 a {\em weak} dynamic
 programming principle  for Problem  \eqref{Value} when $g$ is assumed to be  Borelian only, which follows from some fine results recently obtained in \cite{DQS3} together with some properties of DRBSDEs.  We then focus on the continuous case for which we prove a {\em strong} dynamic programming principle.
 We stress that it cannot be directly derived from the  {\em weak} dynamic programming principle. To this aim, we show  in particular the continuity of the value function with respect to time $t$, which requires some refined properties of doubly reflected BSDEs with jumps. The last part of the paper is devoted to the relation of the value function 
 with generalized Hamilton--Jacobi--Bellman variational inequalities (HJBVIs). 
 We prove that the value function is a {\em weak} viscosity solution of the HJBVIs 
 in the irregular case,  and a {\em classical} viscosity solution in the continuous case. Uniqueness is obtained under additional assumptions.  In terms of PDEs, our result  provides a probabilistic interpretation of nonlinear HJBVIs.
 
 This work completes the one of Buckdahn-Li \cite{Buckdahn2}),
who  have studied  a related optimization problem for doubly reflected BSDEs (of the form 
$\sup_{\alpha} \inf_{\beta}$)
in the case of a Brownian filtration and a continuous reward $g$. Unlike our approach, they do not use a dynamic programming principle to show that the value function  is a viscosity solution of a generalized HJB equation.

 The paper is organized as follows. In Section \ref{sec3}, we introduce the mixed {\em generalized} Dynkin game/stochastic control problem. In Section \ref{sec-prop}, we provide some preliminary properties for doubly reflected BSDEs with jumps. 
 In Section \ref{section4}, we prove the dynamic programming principles both in the irregular and the regular case. 
 In Section \ref{sec5}, we derive that the value function of our problem is a weak (respectively classical) viscosity solution 
 of some generalized HJBVIs in the irregular (respectively regular) case. In Section \ref{comp}, using the results obtained in the previous sections, we provide some additional results for the value function in the discontinuous 
 case.
 In the Appendix, we provide some complementary properties which are used in the paper.
 
\paragraph{Related literature on games and  mixed control problems with stopping times.}

Classical Markovian Dynkin games (with linear expectation) have been studied in particular by Bensoussan-Friedman \cite{BF},  and by 
 Bismut \cite{Bismut}, Alario-Nazaret et al.\cite{ALM}, Kobylanski et al  \cite{KQC} in a non Markovian framework. 
See also  e.g. Hamad\`ene-Lepeltier \cite{HL} for the study of mixed classical Dynkin games. 
Links between Dynkin games and 
doubly reflected BSDEs have been provided in the classical case (see e.g.  Cvitani\'c-Karatzas \cite{CK}, Hamad\`ene and Hassani \cite{HH}, Hamad\`ene and Wang \cite{HW}, Hamad\`ene and Ouknine \cite{HO}),  
and extended to {\em generalized} Dynkin games (that is with nonlinear expectation) in \cite{DQS2}.

 Controller/stopper games are special cases of  Problem \eqref{Value}  when player $A$ is only a controller (no stopping time $\tau$). 
  They  have been studied  in the case of linear expectation by  e.g.  Karatzas-Zamfirescu \cite{KZ}, Bayraktar-Huang \cite{BH} and Choukroun et al. \cite{CCP}.

Finally, mixed optimal stopping/stochastic control problems  are special cases of  Problem \eqref{Value} when there is no player $B$. 
There is   then no more game  aspect. In the particular   case of linear expectation, we refer to  Bensoussan-Lions \cite{BL}, \O ksendal-Sulem
\cite{OS},  
and Bouchard-Touzi \cite{BT} who 
provided  a weak dynamic programming principle when the value function is  irregular. 
This result has been extended to the ${\cal E}^f$-expectation case in  \cite{DQS3}, where  links between the value function and  generalized HJBIs are 
provided under very weak assumptions on the terminal cost (or reward)  function.

We have seen above that Problem \eqref{Value2} is related to an optimization problem on DRBSDEs. 
When the terminal reward map $g$ is continuous, 
some optimization problems relative to BSDEs (see e.g. Peng \cite{Peng92}), to RBSDEs (see Buckdahn-Li \cite{BuckdahnR}) or 
DRBSDEs
(Buckdahn-Li \cite{Buckdahn2} in the Brownian case)
 have been studied in the literature.

\paragraph{Motivating applications in mathematical finance.} 
 Links between classical Dynkin games and game options have been provided by e.g. 
Kifer \cite{Kifer}, Kifer and Yu \cite{KiferYu}, Hamad\`ene \cite{H}.
Recall that a game contingent claim is a contract between a seller and a buyer 
which allows the seller to cancel it at a stopping time time $\sigma \in {\cal T}$ and the buyer to exercise it at any time 
$\tau \in {\cal T}$. The process $X$ may be interpreted as the price process of the underlying asset. 
If the buyer (resp. the seller)  exercises (resp. cancels) at maturity time $T$, then   the seller pays the amount $g(X_{T})$  to the buyer.
 If the buyer exercises at time $\tau <T$ before the seller cancels, then the seller pays the buyer  the amount $ h_1( X_{\tau})$, but if the seller 
cancels  before the buyer exercises, then he pays  the amount $h_2( X_{\sigma})$ to the buyer at the cancellation time $\sigma$. The difference 
$h_2( X_{\sigma}) - h_1( X_{\sigma}) \geq 0$ is interpreted as a penalty that the seller pays to the buyer for the cancellation of the contract. 
In other terms, if the seller (resp. the buyer)  selects a cancellation time $\sigma$ (resp. an exercise time $\tau$),  the seller pays to the buyer at time 
$\tau \wedge \sigma$
the {\em payoff}
$$ I(\tau, \sigma) := h_1( X_{\tau})\textbf{1}_{ \tau \leq \sigma, \tau <T }+
h_2( X_{\sigma})\textbf{1}_{\sigma < \tau }+
g(X_{T})\textbf{1}_{\tau = \sigma=T}].$$
In a perfect market model, there exists an (unique) fair price  for a game
 option, which can be characterized as the value function of a  Dynkin game of the form
 $$\sup_{\tau \in \mathcal{T} }
\inf_{\sigma \in \mathcal{T}} I(\tau, \sigma)= \sup_{\tau \in \mathcal{T} }
\inf_{\sigma \in \mathcal{T}} {\bf E} [ h_1( X_{\tau})\textbf{1}_{ \tau \leq \sigma, \tau <T }+
h_2( X_{\sigma})\textbf{1}_{\sigma < \tau }+
g(X_{T})\textbf{1}_{\tau = \sigma=T}],$$
 where the expectation $ {\bf E} $ is taken under the risk-neutral measure, and 
 $X$ may be interpreted as the price process of the underlying asset (see \cite{Kifer},  \cite{KiferYu} and \cite{H}). Recently, we have generalized in \cite{DQS4} this result  to the case of nonlinear pricing by using  the results on {\em generalized} Dynkin games  
 provided 
 in  \cite{DQS2}. The fair price of the game option is then of the form \eqref{optionjeu}.

 In the presence of constraints and ambiguity on the model represented by an ``ambiguity" parameter $\alpha \in {\mathcal A}$, the controlled payoff  $I^\alpha(\tau,\sigma)$ is given by \eqref{Perf}, and the nonlinear expectation by
 $\mathcal{E}^{\alpha}$. In this case, 
 there is a set of possible fair prices for the game option which can be written as  
 $$\sup_{\tau \in \mathcal{T}}
\inf_{\sigma \in \mathcal{T}}\mathcal{E}_{0, \tau \wedge \sigma}^{\alpha} [ I^\alpha(\tau,\sigma)]= \sup_{\tau \in \mathcal{T}}
\inf_{\sigma \in \mathcal{T}}\mathcal{E}_{0, \tau \wedge \sigma}^{\alpha} [ h_1( X_{\tau}^{\alpha})\textbf{1}_{ \tau \leq \sigma, \tau <T }+
h_2( X_{\sigma}^{\alpha})\textbf{1}_{\sigma < \tau }+
g(X_{T}^{\alpha})\textbf{1}_{\tau = \sigma=T}],$$
where $\alpha \in {\mathcal A}$.
 The supremum of these possible prices  over $ {\mathcal A}$ then corresponds to the value function of 
 the mixed {\em generalized} Dynkin game/stochastic control problem  \eqref{Value2}.
 
\section{Mixed stochastic control/generalized Dynkin game}\label{sec3}

 Let $T>0$ be fixed.  
 We consider the product space $\Omega:=\Omega_W \otimes \Omega_N$, where $\Omega_W:= {\cal C} ([0,T])$ is the  Wiener space, that is the set of continuous functions $\omega^1$ from $[0,T]$ into $\mathbb{R}^p$ such that $\omega^1(0) = 0$, and $\Omega_N := \mathbb{D}([0,T])$ is the Skorohod space of right-continuous with left limits (RCLL)  functions $\omega^2$ from $[0,T]$ into  $\mathbb{R}^d$,  such that 
 $\omega^2(0) = 0$.   
Recall that $\Omega$ is a Polish space for the topology of Skorohod. 
Here  $p, d \geq 1$, but, for notational simplicity, we shall consider only $\mathbb{R}$-valued functions, that is the case $p=d=1$. 

Let $B= (B^1, B^2)$ be the canonical process defined for each $t \in [0,T]$ and each $\omega= (\omega^1, \omega^2)$ by $B^i_t(\omega)= B^i_t(\omega^i):=\omega^i_t$, for $i=1,2$.   
Let us denote the first coordinate process $B^1$ by $W$. 
Let $P^W$ be the probability measure on  $(\Omega_W,\mathcal{B}(\Omega_W))$ such that $W$ is a Brownian motion. Here $\mathcal{B}(\Omega_W)$ denotes the Borelian $\sigma$-algebra on $\Omega_W$.

Set ${\bf E}:= \mathbb{R}^n \backslash \{0\}$ equipped with its Borelian $\sigma$-algebra $\mathcal{B}(\bf{E})$, where $n \geq 1$, and a $\sigma$-finite positive measure $\nu$ such that  $\int_{\bf E} (1 \wedge |e| ) \nu(de) < \infty$.
We  define 
 the jump random measure $N$ as follows: for each $t>0$ and each ${\bf B}$ $ \in \mathcal{B}(\bf{E})$, 
\begin{equation}\label{measure}
N(\omega, [0,t] \times {\bf B})= N(\omega^2, [0,t] \times {\bf B}):= \sum_{0< s \leq t} \textbf{1}_{\{\Delta \omega_s^2 \in {\bf B} \}}.
\end{equation}
Let $P^N$ be the probability measure on  $(\Omega_N,\mathcal{B}(\Omega_N))$ 
 such that $N$ is a Poisson random measure with compensator $\nu(de)dt$ and such that 
 $B_t^2= \sum_{0< s \leq t} \Delta B_s^2 $ a.s. Note that the sum of  jumps  is well defined up to a $P^N$-null set.
 We denote by $\Tilde{N}(dt,de):={N}(dt,de) -  \nu(de)dt$ the compensated Poisson measure.
  The space $\Omega$ is equipped with its Borelian $\sigma$-algebra $\mathcal{B}(\Omega)$ and the probability measure $P:= P^W \otimes P^N$.
   Let $\mathbb{F}:=(\mathcal{F}_t)_{t \leq T}$ be  the filtration generated by $W$ and $N$ {\em completed} 
with respect 
to $\mathcal{B}(\Omega)$ and $P$ (see e.g. \cite{J} p.3 for the definition of a {\em completed} filtration).  Note that ${\cal F}_T$ is equal to the {\em completion} of the $\sigma$-algebra $\mathcal{B}(\Omega)$ with respect to $P$, and 
${\cal F}_0$ is the $\sigma$-algebra generated by $P$-null sets.
Let ${\cal P}$ be the predictable $\sigma$-algebra on $\Omega \times [0,T]$ associated 
with the filtration $\mathbb{F}$.

%
%
%

  We introduce the following spaces : \begin{itemize}
\item $\mathbb{H}^{2}_T$ (also denoted  $\mathbb{H}^{2}$)  := the set of real-valued 
predictable processes ($Z_t$) with $\mathbb{E}\int_0^T Z_s^{2}ds < \infty$. 
\item $\mathcal{S}^2$ :=  the set of real-valued RCLL adapted processes $(\varphi_s)$  with $\mathbb{E}[\sup_{0\leq s \leq T}\varphi_s^2] < \infty. $ 
\item  ${L}_{\nu}^2$ := the set of measurable functions $l:({\bf E},\mathcal{K}) \rightarrow (\mathbb{R}, {\cal B} (\mathbb{R}))$  such that \\
$\| l \|_\nu^2 :=\int_{{\bf E}}l^2(e)\nu(de)< \infty.$
The set  ${L}^2_\nu$ is a 
 Hilbert space equipped with the scalar product 
$\langle l    , \, l' \rangle_\nu := \int_{{\bf E}} l(e) l'(e) \nu(de)$ for all $ l , \, l' \in {L}^2_\nu \times {L}^2_\nu.$ 
\item
$\mathbb{H}^{2}_{\nu}$ : the set of predictable real-valued processes $(k_t(\cdot))$ with $\mathbb{E}\int_0^T \|k_s\|_{{L}_{\nu}^2}^2ds < \infty$. 
\end{itemize}

Let ${\bf A}$ be a nonempty closed subset of ${\mathbb R}^p$.
 Let   $\mathcal{A}$ be the set of controls,  defined as the set of predictable processes $\alpha$  valued in 
 ${\bf A}$.
For each $\alpha \in \mathcal{A}$, initial time  $t \in [0,T]$ and initial condition $x$ in $\mathbb{R}$,  let $( X_s^{\alpha,t,x})_{ t \leq s \leq T}$ be the unique $\mathbb{R}$-valued solution in $\mathcal{S}^2$ of the stochastic differential equation:
\begin{equation}\label{richesse}
X_s^{\alpha,t,x}=\dd x+ \int_t^s b(X_r^{\alpha,t,x}, \alpha_r)dr+\int_t^s \sigma(X_r^{\alpha,t,x}, \alpha_r)dW_r+\int_t^s \int_{{\bf E}} \beta(X_{r^{-}}^{\alpha,t,x}, \alpha_r,e) \Tilde{N}(dr,de),
\end{equation}
where  $b, \ \sigma :\mathbb{R} \times {\bf A} \rightarrow \mathbb{R}$, 
are  Lipschitz continuous with respect to $x$ and $\alpha$,   and $\beta : \mathbb{R} \times {\bf A} \times {\bf E} \rightarrow \mathbb{R}$ is  a bounded measurable function such that for some constant $C \geq 0$, and for all $e \in \mathbb{R}$
\begin{align*}
&|\beta(x,\alpha,e)| \leq C\,\Psi(e),
  \;\; x \in \mathbb{R}, \alpha \in {\bf A}\\
&|\beta(x,\alpha,e)- \beta(x',\alpha',e)| \leq C(|x-x'| + |\alpha - \alpha'|)
\Psi(e), 
\;\;  x, x'\in \mathbb{R}, \alpha, \alpha' \in {\bf A},
\end{align*}
where $\Psi \in $  ${L}^2_\nu \cap {L}^1_\nu$. 


 The criterion of our mixed control problem, depending on $\alpha$, is defined via  a BSDE with a Lipschitz driver function  $f$ satisfying the following conditions:

\begin{itemize}
\item[(i)]
$f: {\bf A} \times [0,T] \times  \mathbb{R}^3 \times  {L}_{\nu}^2  \rightarrow (\mathbb{R}, \mathcal{B}(\mathbb{R}))$ is  $ 
\mathcal{B}({\bf A}) \otimes \mathcal{B}([0,T]) \otimes \mathcal{B}(\mathbb{R}^3) \otimes \mathcal{B}({L}_{\nu}^2)$-measurable
\item[(ii)]
 $|f(\alpha,t,x,0,0,0)| \leq C(1+|x|^p), \forall \alpha \in {\bf A}, t \in [0,T],   x\in\mathbb{R}$, where  $p \in 
 \mathbb{N}^*$.

\item[(iii)]
$ |f(\alpha,t,x,y,z,k)- f (\alpha',t,x',y',z',k')| \leq C(  |\alpha -\alpha '|+ |x-x'|+|y-y'|+|z-z'|+\|k-k'\|_{{L}_{\nu}^2})$, $\forall t \in [0,T]$, $x, x', y, y', z,z' \in \mathbb{R}$, $k,k' \in {L}_{\nu}^2, \alpha, \alpha' \in {\bf A}.$

\item[(iv)]
$ f(\alpha,t,x,y,z,k_2)- f (\alpha,t,x,y,z,k_1) \geq <\gamma(\alpha,t,x,y,z,k_1,k_2),k_2-k_1>_{\nu}, \forall t,x,y,z,k_1,k_2, \alpha,$
\end{itemize}
where $\gamma: {\bf A} \times [0,T] \times  \mathbb{R}^3 \times  ({L}_{\nu}^2)^2  \rightarrow {L}_{\nu}^2$ is $ 
\mathcal{B}({\bf A}) \otimes \mathcal{B}([0,T]) \otimes \mathcal{B}(\mathbb{R}^3) \otimes \mathcal{B}(({L}_{\nu}^2)^2)/ 
 \mathcal{B}({L}_{\nu}^2)$-measurable,  
$  |\gamma(.)(e)| \leq  \Psi(e) 
 $ and $\gamma(.)(e)\geq -1$ $d\nu(e)$-a.s. where $\Psi \in {L}^2_\nu$.

Condition (iv) allows us to apply the comparison theorem for non reflected BSDEs, reflected and doubly reflected BSDEs with jumps (see \cite{16}, \cite{17} and \cite{DQS2}).

Let $(t,x) \in [0,T] \times \mathbb{R}$ and $\alpha \in \mathcal{A}$. Let $\mathcal{E}^{f^{\alpha,t,x}}$ 
(also denoted  by $\mathcal{E}^{\alpha, t, x}$)  be the nonlinear 
conditional expectation  associated with $f^{\alpha,t,x}$, 
defined for each stopping time $S$ and for each  $\zeta \in {L}^2(\mathcal{F}_S)$ as:
$$\mathcal{E}^{\alpha,t,x}_{r,S}[\zeta]:=y_r^{\alpha,t,x},  \;\; t \leq r \leq S,$$
where $(y_r^{\alpha,t,x})_{t \leq r \leq S}$ is the solution in $\mathcal{S}^2$ of the BSDE associated with driver $f^{\alpha,t,x}(r, y,z,k):= f(\alpha_r,r, X_r^{\alpha,t,x},y,z,k)$, terminal time $S$ and terminal condition $\zeta$, that is satisfying the dynamics
\begin{equation}\label{for}
-dy_r^{\alpha,t,x}=f(\alpha_r,r,X_r^{\alpha,t,x},y_r^{\alpha,t,x},z_r^{\alpha,t,x},k_r^{\alpha,t,x})dr-z_r^{\alpha,t,x}dW_r-\int_{{\bf E}}k_r^{\alpha,t,x}(e)\Tilde{N}(dr,de),
\end{equation}
with $y_S^{\alpha,t,x}=\zeta$, 
where $z_\cdot^{\alpha,t,x}, k_\cdot^{\alpha,t,x}$ 
are the associated processes in $\mathbb{H}^{2}$ and  $\mathbb{H}^{2}_{\nu}$ respectively.
%

 For all $(t,x) \in [0,T] \times \mathbb{R}$ and all control $\alpha \in \mathcal{A}$, 
 we define the barriers for $i={1,2}$ by
$h_{i}(s,X_s^{\alpha,t,x}),$ for $t \leq s < T$,  and the terminal condition by $g(  X_T^{\alpha,t,x})$, where 
%
\begin{itemize}
\item[(i)]
$g :  \mathbb{R} \rightarrow \mathbb{R}$ is Borelian,

%
%
%
%
%
%
\item[(ii)]
 $h_1 :[0,T] \times \mathbb{R} \rightarrow \mathbb{R}$ and 
 $h_2 :[0,T] \times \mathbb{R} \rightarrow \mathbb{R}$ are
   functions which are Lipschitz continuous with respect to $x$ uniformly in $t$, and 
   continuous with respect to $t$ on $[0,T]$, with $h_1 \leq h_2$,
   \item[(iii)] $h_1$ (or $h_2$) is ${\cal C}^{1,2}$ with bounded derivatives,
\item[(iv)]
$\quad    |h_1(t,x)|+|h_2(t,x)| +  |g(x)| \leq C(1+ |x|^p), \forall t \in [0,T], x \in \mathbb{R},  \text{ with $p\in \mathbb{N}$}. $
\end{itemize}

 Let $\mathcal{T}$ be the set of stopping times with values in $[0,T]$. 
Suppose the initial time is equal to $0$. For each initial condition $x \in \mathbb{R}$, we consider 
the following  mixed {\em generalized} Dynkin game and stochastic control problem:
\begin{equation}\label{probform}
u(0,x):= \sup_{\alpha \in \mathcal{A}} \sup_{\tau \in \mathcal{T}} 
\inf_{\sigma \in \mathcal{T}} \mathcal{E}_{0, \tau \wedge \sigma}^{\alpha,0,x} 
\left[ h_1(\tau, X_{\tau}^{\alpha,0,x})\textbf{1}_{ \tau \leq \sigma, \tau <T }+
h_2(\sigma, X_{\sigma}^{\alpha,0,x})\textbf{1}_{\sigma < \tau }+
g(X_{T}^{\alpha,0,x})\textbf{1}_{\tau \wedge \sigma=T}\right].
\end{equation}

We now make the problem dynamic.  We  define, for $t \in [0,T]$ and each $\omega$ $\in$ $\Omega$ the $t$-translated path  $\omega^t = (\omega^t_s)_{s \geq t}:=(\omega_s-\omega_t)_{s \geq t}$.
Note that $(\omega^{1,t}_s)_{s \geq t}:=(\omega_s^1-\omega_t^1)_{s \geq t}$  corresponds to the realizations of the translated Brownian motion $W^t:= (W_s-W_t)_{s \geq t}$ and that the translated Poisson random measure $N ^t :=
N(]t,s],.)_{s \geq t}$ can be expressed in terms of $(\omega^{2,t}_s)_{s \geq t}:=(\omega_s^2-\omega_t^2)_{s\geq t}$
 similarly to \eqref{measure}.  Let $\mathbb{F}^t=(\mathcal{F}_s^t)_{t \leq s \leq T}$ be the  filtration 
generated by $W^t$ and $N ^t$ {\em completed with respect to ${\cal B}(\Omega)$ and $P$}. 
Note that for each $s\in [t,T]$, $\mathcal{F}_s^t$ is the $\sigma$-algebra generated by $W_r^t$, $N_r ^t$, $t\leq r\leq s$  and $\mathcal{F}_0$. Recall also that  we have a martingale representation theorem for 
  ${\mathbb F}^t$-martingales as stochastic integrals with respect to $ W^t$ and $\Tilde{N}^t$.

Let us denote by $\mathcal{T}^t_{t}$ the set of stopping times with respect to $\mathbb{F}^t$ with values in $[t,T]$.  Let   $\mathcal{P}^t$ be the predictable $\sigma$-algebra on 
$\Omega \times [t,T]$ equipped with the filtration $\mathbb{F}^t$. \\
We introduce the following spaces of processes. Let $t \in [0,T]$. \\
 Let $\mathbb{H}_t^2$ be  the ${\cal P}^t$-measurable processes $Z$  on $\Omega \times [t,T]$ such that 
 $\| Z \|_{\mathbb{H}_t^2} := \mathbb{E}[\int_t^T Z_u^2du] < \infty.$ \\
 Let 
 $\mathbb{H}_{t, \nu}^2$ be  the set of ${\cal P}^t$-measurable processes $K$  on $\Omega \times [t,T]$ 
with 
 $\|K \|_{\mathbb{H}_{t, \nu}^2} := \mathbb{E}[\int_t^T  \|K_u\|_{\nu}^2du] < \infty.$ \\
 Let  ${\cal S}^2_t$ be  the set of ${\mathbb R}$-valued RCLL processes $\varphi$ on $\Omega \times [t,T]$, ${\mathbb F}^t$-adapted, with $\mathbb{E}[\sup_{t\leq s \leq T} \varphi_s^2] < \infty. $

  Let   $\mathcal{A}^t_t$ be the set of controls $\alpha:\Omega \times [t, T] \mapsto {\bf A}$, which 
 are ${\cal P}^t$-measurable. 
For each initial time $t$ and each initial condition $x$, the value function is defined by:
\begin{equation}\label{probform1}
u(t,x):= \sup_{\alpha \in \mathcal{A}^t_t} \sup_{\tau \in \mathcal{T}^t_{t}} \inf_{\sigma \in  \mathcal{T}^t_{t}} \mathcal{E}_{t, \tau \wedge \sigma}^{\alpha,t,x}\left[ h_1(\tau, X_{\tau}^{\alpha,t,x})\textbf{1}_{\tau
 \leq  \sigma, \tau <T}+h_2(\sigma, X_{\sigma}^{\alpha,t,x})\textbf{1}_{\sigma < \tau}+g(X_{T}^{\alpha,t,x})
 \textbf{1}_{\tau \wedge \sigma=T}\right].  
\end{equation}
Note that since $\alpha$, $\tau$ and $\sigma$ depend only on $\omega^t$, the SDE \eqref{richesse} and the BSDE \eqref{for}
 can be solved with respect to the translated Brownian motion $(W_s-W_t)_{s \geq t}$ and the translated Poisson random measure $N(]t,s], \cdot)_{s \geq t}$. 
Hence the function $u$ is well defined as a deterministic function of $t$ and $x$. 

For each $\alpha \in  \mathcal{A}^t_t$, we introduce the function $u^{\alpha}$ defined as 
 $$u^{\alpha}(t,x):= \sup_{\tau \in \mathcal{T}^t_{t}} \inf_{\sigma \in  \mathcal{T}^t_{t}} \mathcal{E}_{t, \tau \wedge \sigma}^{\alpha,t,x}\left[ h_1(\tau, X_{\tau}^{\alpha,t,x})\textbf{1}_{\tau
 \leq  \sigma, \tau <T}+h_2(\sigma, X_{\sigma}^{\alpha,t,x})\textbf{1}_{ \sigma < \tau}+g(X_{T}^{\alpha,t,x})\textbf{1}_{\tau \wedge \sigma= T}\right].$$
 We thus get 
 \begin{equation}\label{eg}
 u(t,x) =\sup_{\alpha \in \mathcal{A}^t_t} u^{\alpha}(t,x).
 \end{equation}
 Note that for all $\alpha,x$ and $t<T$, we have $h_1(t,x) \leq u^{\alpha}(t,x) \leq h_2(t,x)$,  and hence $h_1(t,x) \leq u(t,x) 
 \leq h_2(t,x)$. 
 Moreover, $u^{\alpha}(T,x) = u(T,x)= g(x)$. 
 
  By Assumption $(iii)$, $h_1$ (or $h_2$) is ${\cal C}^{1,2}$ with bounded derivatives. It follows 
  from Proposition \ref{Moko} in the Appendix, that for each $t \in [0,T]$ and for each $\alpha$ $\in$ ${\cal A}_t^t$, the processes  $\xi_s^{\alpha,t,x} :=h_1(s,X_s^{\alpha,t,x}){\bf 1}_{s<T} + g(X^{\alpha,t,x}_T){\bf 1}_{s=T}$ and $\zeta_s ^{\alpha,t,x}:=h_2(s,X_s^{\alpha,t,x}){\bf 1}_{s<T} + g(X^{\alpha,t,x}_T){\bf 1}_{s=T}$ 
 satisfy Mokobodzki's condition: 
there exist  two  nonnegative  ${\mathbb F}^t$-supermartingales $H^{\alpha,t,x}$ and $H^{' \alpha,t,x}$ in ${\cal S}_t^2$ such that 
 \begin{equation}\label{mo}
 \xi_s^{\alpha,t,x} \leq  H^{\alpha,t,x}_s - H^{' \alpha,t,x}_s \leq \zeta_s ^{\alpha,t,x}, 
 \,\,\,t \leq s \leq T \quad {\rm a.s.}
 \end{equation}
%
%
%
%
By Theorem  4.7 in  \cite{DQS2}, for each $\alpha$, the value function $u^\alpha$ of the above {\em generalized} Dynkin game is characterized as the solution of the doubly reflected BSDE associated with driver $f^{\alpha,t,x}$, 
barriers $\xi_s^{\alpha,t,x}=h_1(s,X_s^{\alpha,t,x})$, $\zeta_s ^{\alpha,t,x}=h_2(s,X_s^{\alpha,t,x})$ for $s<T$, and terminal condition $g(X^{\alpha,t,x}_T)$, that is 
\begin{equation}\label{caractb}
u^\alpha(t,x)=Y_t^{\alpha,t,x},
\end{equation}
where $(Y^{\alpha,t,x},  Z^{\alpha,t,x},   K^{\alpha,t,x})  \in \mathcal{S}^2 \times \mathbb{H}^2 \times
\mathbb{H}^2_{\nu}$ is the solution of the doubly reflected BSDE:
\begin{equation}\label{4.4}
\begin{cases}
Y^{\alpha,t,x}_s=g(X^{\alpha,t,x}_T)+\dd\int_{s}^{T}f(\alpha_r,r,X^{\alpha,t,x}_r,Y^{\alpha,t,x}_r,Z^{\alpha,t,x}_r, K^{\alpha,t,x}_r(\cdot))dr \\
\qquad\dd + A^{1,\alpha,t,x}_T-A^{1,\alpha,t,x}_s- A^{2,\alpha,t,x}_T+A^{2,\alpha,t,x}_s 
-\int_s^T Z^{\alpha,t,x}_rdW_r-\int_s^T \dd\int_{{\bf E}}K^{\alpha,t,x}(r,e) \Tilde{N}(dr,de)\\
\xi^{\alpha,t,x}_s \leq Y^{\alpha,t,x}_s \leq \zeta^{\alpha,t,x}_s, t \leq s < T\; \text{ a.s. },\\
A^{1,\alpha,t,x}, A^{2,\alpha,t,x} \text{ are  RCLL nondecreasing predictable processes with } A^{1,\alpha,t,x}_t =  A^{2,\alpha,t,x}_t=0 \text{ and  }\\
\dd\int_{t}^{T}(Y_{s}^{\alpha,t,x}-\xi^{\alpha,t,x}_{s})dA_s^{1,\alpha,t,x,c}=0 \text{ a.s.   and   } \;  \Delta A^{1,\alpha,t,x,d}_{s}= - \Delta A^{1,\alpha,t,x}_s \,{\bf 1}_{\{Y^{\alpha,t,x}_{s^-} = \xi^{\alpha,t,x}_{s^-}\}}  \quad   \rm{a.s.} \\
\dd\int_{t}^{T}(\zeta_{s}^{\alpha,t,x}-Y_{s}^{2,\alpha,t,x})dA^{2,\alpha,t,x,c}_s=0 \text{ a.s. and  }  \;  \Delta A^{2,\alpha,t,x,d}_{s}= - \Delta A^{2,\alpha,t,x}_s \,{\bf 1}_{\{Y^{\alpha,t,x}_{s^-} = \zeta^{\alpha,t,x}_{s^-}\}}  \quad   \rm{a.s.} \\
\end{cases}
\end{equation}
Here  $A^{1,\alpha,t,x,c}$ (resp. $A^{2,\alpha,t,x,c}$ ) denotes the continuous part of $A^1$ (resp. $A^2$) and $A^{1,\alpha,t,x,d}$ (resp.$A^{2,\alpha,t,x,d}$ )  the discontinuous part.
In the particular case when $h_1(T,x) \leq g(x) \leq h_2(T,x)$, then the obstacles $\xi^{\alpha,t,x}$ and $\zeta^{\alpha,t,x}$ satisfy for all predictable stopping time $\tau$, $\xi_{\tau^-} \leq \xi_{\tau}$ and $\zeta_{\tau^-} \geq \zeta_{\tau}$ a.s.\, 
which implies the continuity 
of $A^{1,\alpha,t,x}$ and $A^{2,\alpha,t,x}$ (see \cite{DQS2}).

Note that the  doubly reflected BSDE \eqref{4.4} can be solved in $\mathcal{S}^2 \times \mathbb{H}_t^2 \times \mathbb{H}_{t, \nu}^2$  with respect to the $t$-translated Brownian motion and the $t$-translated Poisson random measure. 

In the following, for each $\alpha \in \mathcal{A}_t^t$,  $Y^{\alpha,t,x}_s$ will be also denoted by  ${Y}_{s,T}^{\alpha,t,x}[g(X^{\alpha,t,x}_T)]$.


Using \eqref{caractb}, our initial optimization problem \eqref{probform1} can thus be reduced to an optimal  control problem for doubly reflected BSDEs: $$u(t,x) =\sup_{\alpha \in \mathcal{A}^t_t} Y_t^{\alpha,t,x}= \sup_{\alpha \in \mathcal{A}^t_t} {Y}_{t,T}^{\alpha,t,x}[g(X^{\alpha,t,x}_T)].$$

We now provide some new results on  doubly reflected BSDEs, which will be used to prove the dynamic programming principles.

\section{Preliminary properties  for doubly reflected BSDEs}\label{sec-prop}

We show  in a general non Markovian framework a continuity property and a Fatou lemma  for doubly reflected  BSDEs, where the limit involves both terminal condition and terminal time. 

A function $f$ is said to be a {\em Lipschitz driver} if \\
$f: [0,T]  \times \Omega \times \R^2 \times L^2_\nu \rightarrow \R $
$(\omega, t,y, z, k(\cdot)) \mapsto  f(\omega, t,y, z, k(\cdot))  $
  is $ {\cal P} \otimes {\cal B}(\R^2)  \otimes {\cal B}(L^2_\nu) 
- $ measurable,  uniformly Lipschitzian  with respect to $y, z, k(\cdot)$ and such that $f(.,0,0,0) \in \H^2$.

A Lipschitz driver $f$ is said to satisfy Assumption~\ref{Royer} if the following holds:
\begin{assumption}\label{Royer} 
 $dP \otimes dt$-a.s\, for each $(y,z, k_1,k_2)$ $\in$ $ \mathbb{R}^2 \times (L^2_{\nu})^2$,
$$f( t,y,z, k_1)- f(t,y,z, k_2) \geq \langle \gamma_t^{y,z, k_1,k_2}  \,,\,k_1 - k_2 \rangle_\nu,$$ 
with
$
\gamma:  [0,T]  \times \Omega\times \mathbb{R}^2 \times  (L^2_{\nu})^2  \rightarrow  L^2_{\nu}\,; \, (\omega, t, y,z, k_1, k_2) \mapsto 
\gamma_t^{y,z,k_1,k_2}(\omega,.)
$, supposed to be
 ${\cal P } \otimes {\cal B}({\mathbb R}^2) \otimes  {\cal B}( (L^2_{\nu})^2 )$-measurable, and satisfying $ dP\otimes dt $-a.s.\,, for each $(y,z, k_1, k_2)$ $\in$ ${\mathbb R}^2 \times (L^2_{\nu})^2$,
$
\gamma_t^{y,z, k_1, k_2} (\cdot)\geq -1$  and $\|\gamma_t^{y,z, k_1, k_2}(\cdot)\|_{\nu}\leq K$, where $K$ is a positive constant.
\end{assumption}
This assumption ensures  the comparison theorem for BSDEs, reflected BSDEs and doubly reflected BSDEs with jumps  (see \cite{16}, \cite{17} and \cite{DQS2}). It is satisfied if, for example, $f$ is of class ${\cal C}^1$ with respect to $k$ such that 
$ \nabla_k f $ is bounded (in $L^2_{\nu}$) and $ \nabla_k f \geq -1$ (see \cite{DQS2} Proposition A.2).

We extend the definition of reflected and doubly reflected BSDEs when the terminal time 
  is  a stopping time $ \theta$ $ \in {\cal T}$ and the terminal condition is a random variable  $\xi$ in $ L^2({\cal F}_ \theta)$. 
  let $f$ be a given Lipschitz driver. 
Let $(  \eta_t)$ be a given obstacle  RCLL process in $\mathcal{S}^2$. 
%
The solution, denoted $(Y_{., \theta}(\xi), Z_{., \theta}(\xi), k_{., \theta}(\xi))$, of the reflected BSDEs associated with terminal time $\theta$, driver $f$, obstacle $(  \eta_s)_{ s < \theta}$, and terminal condition $\xi$  
 is defined  as the unique solution in $\mathcal{S}^2 \times \mathbb{H}^2 \times
\mathbb{H}^2_{\nu}$ of the reflected BSDE with  terminal time $T$, driver $f(t,y, z,k) {\bf 1}_ { \{ t \leq  \theta  \}  }$, 
 terminal condition $\xi$ and obstacle $ \eta_t {\bf 1}_{t < \theta} + \xi {\bf 1}_{t \geq \theta}$.
 Note that 
$Y_{t, \theta}(\xi) = \xi, Z_{t, \theta}(\xi) =0,k_{t, \theta}(\xi)=0 $ for $t \geq  \theta$.
Similarly, 
let $(  \eta_t)$ and $(\zeta_t)$ be given   RCLL processes in $\mathcal{S}^2$ satisfying  
Mokobodzki 's condition, that is,
there exist  two  nonnegative  supermartingales $H$ and $H^{' }$ in ${\cal S}^2$ such that 
 \begin{equation*}
 \eta_s \leq H_s - H^{'}_s \leq \zeta_s, 
 \,\,\,0 \leq s \leq T \quad {\rm a.s.}
 \end{equation*}
The solution, denoted $(Y_{., \theta}(\xi), Z_{., \theta}(\xi), k_{., \theta}(\xi))$  of the 
doubly reflected BSDEs  associated with terminal stopping time $\theta$, driver $f$, barriers $(  \eta_s)_{ s < \theta}$ and $(  \zeta_s)_{ s < \theta}$, and terminal condition $\xi$,  
 is defined 
 as the unique solution in $\mathcal{S}^2 \times \mathbb{H}^2 \times
\mathbb{H}^2_{\nu}$ of the doubly reflected BSDE with  driver $f(t,y, z,k) {\bf 1}_ { \{ t \leq  \theta  \}  }$,
 terminal time $T$, 
 terminal condition $\xi$ and barriers $ \eta_t {\bf 1}_{t < \theta} + \xi {\bf 1}_{t \geq \theta}$ and $ \zeta_t {\bf 1}_{t < \theta} + \xi {\bf 1}_{t \geq \theta}$.
 Note that 
$Y_{t, \theta}(\xi) = \xi, Z_{t, \theta}(\xi) =0,k_{t, \theta}(\xi)=0 $ for $t \geq  \theta$.

\

We first prove a continuity property for doubly reflected BSDEs where the limit involves both terminal condition and terminal time.

\begin{proposition}\label{continuity}[A continuity property for doubly reflected BSDEs]
Let $T>0$.  Let $f$ be a Lipschitz driver satisfying Assumption \ref{Royer}. Let $(\eta_t)$, $(\zeta_t)$ be two RCLL processes in $\mathcal{S}^2$, with $\eta_. \leq \zeta_.$. Let $f$ be a given Lipschitz driver. 
Let $(\theta^n)_{ n \in \mathbb{N}} $ be a non increasing sequence of 
stopping times in 
$\mathcal{T}$,  converging a.s. to  $\theta \in \mathcal{T}$ as 
$n$ tends to $\infty$. Let   $(\xi^{n})_{ n \in \mathbb{N}} $ be a sequence of random variables such that 
$E[{\rm ess }\sup_{n}( \xi^n)^2]< + \infty$, and  for each $n$,
$\xi^{n}$ is ${\mathcal F}_{\theta^{n}}$-measurable.  Suppose that for each $n$, the processes $\eta_t {\bf 1}_{\{t< \theta^n\}} + \xi^n {\bf 1}_{\{t\geq \theta^n\}}$ and $\zeta_t {\bf 1}_{\{t< \theta^n\}} + \xi^n {\bf 1}_{\{t\geq \theta^n\}}$ satisfy Mokobodzki 's condition.
Suppose  that $\xi^{n}$ converges a.s. to an ${\mathcal F}_{\theta}$-measurable random variable $\xi$ as 
$n$ tends to $\infty$.
Suppose that 
\begin{equation}\label{ass}
\eta_{\theta} \leq \xi \leq \zeta _{\theta} \quad {\rm a.s.}\,,
\end{equation}
and that the processes $\eta_t {\bf 1}_{\{t< \theta\}} + \xi {\bf 1}_{\{t\geq \theta\}}$ and $\zeta_t {\bf 1}_{\{t< \theta\}} + \xi {\bf 1}_{\{t\geq \theta\}}$ satisfy Mokobodzki 's condition.

Let $Y_{.,\theta^{n}}(\xi^n )$; $Y_{.,\theta}(\xi )$ be the solutions of the doubly reflected BSDEs associated with driver $f$, barriers $(\eta_s)_{ s < \theta^n}$ and $(\zeta_s)_{ s < \theta^n}$(resp. $(\eta_s)_{ s < \theta}$ and $(\zeta_s)_{ s < \theta}$) , terminal time $\theta^n$ (resp. $\theta$), terminal condition $\xi^n$ (resp. $\xi$).
Then,  for each  stopping time $\tau$ with $\tau \leq \theta$ a.s.\,,  
 $$Y_{\tau,\theta}(\xi ) = \lim_{n \rightarrow + \infty} Y_{\tau,\theta^{n}}(\xi^n ) \quad a.s.$$
When for each $n$,  $\theta_n= \theta$ a.s.\,, the result still holds without Assumption \eqref{ass}.

\begin{remark}
 As in the case of reflected BSDEs (see Proposition \ref{CRBSDE} in the Appendix),  there is an extra difficulty due to the presence 
 of the barriers (and the variation of the terminal time). The additional assumption  \eqref{ass} on the obstacle is here required to obtain the result. 
\end{remark}

 \end{proposition}

  \dproof  We first consider the simpler case when for each $n$, $\theta_n = \theta$ a.s. 
  By the a priori estimates on doubly reflected BSDEs provided in  \cite{DQS2} (see Proposition 5.3) and
  the convergence of $\xi_n$ to $\xi$, we derive that  $Y_{\tau,\theta}(\xi ) = \lim_{n \rightarrow + \infty} Y_{\tau,\theta}(\xi^n ) $ a.s.\\ 
 We now turn to the general case. In this  case, Proposition 5.3 in \cite{DQS2} does not give the result. 
 By the flow property  for doubly reflected BSDEs (or ``semigroup property", see \cite{Buckdahn2}), we have:
  $$
 Y_{\tau, \theta_n}(\xi^n)=Y_{\tau, \theta}(Y_{\theta, \theta_n}(\xi^n)).
 $$
 By the first step,  it is thus sufficient to show that $\lim_{n \rightarrow} Y_{\theta, \theta_n}(\xi^n)=\xi$ a.s.\\
 Since the solution of the doubly reflected BSDE associated with terminal condition  $\xi^n$ and terminal time $\theta_n$ is smaller than the solution $\tilde Y_{., \theta_n}(\xi^n)$ of the reflected BSDE associated with (one) obstacle $(\xi_t)_{t < \theta_n}$, terminal condition  $\xi^n$ and terminal time $\theta_n$ , we have:
 \begin{equation}\label{ineq}
 Y_{\theta, \theta_n}(\xi^n) \leq \tilde Y_{\theta, \theta_n}(\xi^n)  \text{ a.s. }
 \end{equation}
 By the continuity property of reflected BSDEs with respect to terminal time and terminal condition 
 (Lemma \ref{CRBSDE}), we have $\lim_{n \rightarrow \infty} \tilde Y_{\theta, \theta_n}(\xi^n)= \tilde Y_{\theta, \theta}(\xi)= \xi$ a.s. Taking the $\lim \sup$ in $\eqref{ineq}$, we get
 $$
 \lim \sup_{n \rightarrow \infty} Y_{\theta, \theta_n}(\xi^n) \leq \lim_{n \rightarrow \infty} \tilde Y_{\theta, \theta_n}(\xi^n) = \xi \,\,\, \text{ a.s.}
$$
  It remains to show that 
   \begin{equation}\label{dr}
 \lim \inf_{n \rightarrow \infty} Y_{\theta, \theta_n}(\xi^n) \geq  \xi \,\,\, \text{ a.s.}
 \end{equation}
Since the solution of the doubly reflected BSDE associated with terminal condition  $\xi^n$ and terminal time $\theta_n$ is greater  than the solution $\bar Y_{., \theta_n}(\xi^n)$ of the reflected BSDE associated with the upper obstacle $(\zeta_t)_{t < \theta_n}$, terminal condition  $\xi^n$ and terminal time $\theta_n$, we have:
 \begin{equation}\label{ineq2}
 Y_{\theta, \theta_n}(\xi^n) \geq \bar Y_{\theta, \theta_n}(\xi^n)  \text{ a.s. }
 \end{equation}
 By Lemma \ref{CRBSDE}, we derive that the solution of a reflected BSDEs with upper barrier is continuous with respect to both terminal time and terminal condition.  Hence, $\lim_{n \rightarrow \infty} \bar Y_{\theta, \theta_n}(\xi^n)= \xi$ a.s. Taking the $\lim \inf$ in $\eqref{ineq2}$, we derive 
 inequality \eqref{dr}.
  \fproof
 
%

Using Proposition \ref{continuity} together with the monotonicity property of the solution of a doubly reflected BSDEs with respect to terminal condition, one can derive the following Fatou lemma.
\begin{lemma}[A Fatou lemma for doubly reflected BSDEs]\label{fatou2}
Let $T>0$.  Let $(\eta_t)$, $(\zeta_t)$ be two RCLL processes in $\mathcal{S}^2$, with $\eta_. \leq \zeta_.$, and satisfying Mokobodzki's condition. Let $f$ be a Lipschitz driver satisfying Assumption \ref{Royer}. 
Let $(\theta^n)_{ n \in \mathbb{N}} $ be a non increasing sequence of stopping times in 
$\mathcal{T}$,  converging a.s. to  $\theta \in \mathcal{T}$ as 
$n$ tends to $\infty$. Let  $(\xi^{n})_{ n \in \mathbb{N}} $ be a sequence of random variables such that 
$E[\sup_{n}( \xi^n)^2]< + \infty$, and  for each $n$,
$\xi^{n}$ is ${\mathcal F}_{\theta^{n}}$-measurable.
Suppose that for each $n$, the processes $\eta_t {\bf 1}_{\{t< \theta^n\}} + \xi^n {\bf 1}_{\{t\geq \theta^n\}}$ and $\zeta_t {\bf 1}_{\{t< \theta^n\}} + \xi^n {\bf 1}_{\{t \geq \theta^n\}}$ satisfy Mokobodzki's condition.

Let $\underline{\xi}:= \liminf_{n \rightarrow + \infty} \xi^{n}$. Suppose that the processes $\eta_t {\bf 1}_{\{t< \theta\}} + \underline{\xi} {\bf 1}_{\{t\geq  \theta\}}$ and $\zeta_t {\bf 1}_{\{t< \theta\}} + \underline{\xi} {\bf 1}_{\{t\geq  \theta\}}$ satisfy Mokobodzki's condition, and  that \begin{equation}\label{ij}
\eta_{\theta} \leq \, \underline{\xi}\, \leq \zeta_{\theta} \quad {\rm a.s.}
\end{equation}

Let $Y_{.,\theta^{n}}(\xi^n )$ (resp. $Y_{.,\theta}(\, \underline{\xi}\,)$) be the solution(s) of the doubly reflected BSDEs associated with driver $f$, barriers $(\eta_s)_{ s < \theta^n}$ and $(\zeta_s)_{ s < \theta^n}$  (resp. $(\eta_s)_{ s < \theta}$ and $(\zeta_s)_{ s < \theta}$), terminal time $\theta^n$ (resp. $\theta$), terminal condition $\xi^n$ (resp.  $\underline{\xi}$).

Then,  for each  stopping time $\tau $ with $\tau \leq \theta$ a.s.\,,  
\begin{equation}\label{inegalite}
Y_{\tau,\theta}(\, \underline{\xi}\, ) \leq \liminf_{n \rightarrow + \infty} Y_{\tau,\theta^{n}}(\xi^n ) 
  \quad {\rm a.s.}\,
 \end{equation}
When for each $n$,  $\theta_n= \theta$ a.s.\,, the result still holds without Assumption \eqref{ij}.
  \end{lemma}
 \begin{remark}\label{limsup}
 The result also holds for $\bar{\xi}:=\limsup_{n \rightarrow + \infty} \xi^{n}$ instead of $\underline{\xi}$
 with the converse inequality
  in \eqref{inegalite}, \\
  {\bf ICI liminf replaced by limsup,}\\
  under the same assumptions with $\underline{\xi}$ replaced by $\bar{\xi}$ 
 \end{remark}
\noindent The proof, which is very similar to that of the Fatou lemma for classical BSDEs (see Lemma \ref{fatou1}), is left to the reader.
This  lemma  will be used to obtain a weak dynamic principle.
 
\section{Dynamic programming principles}\label{section4}
We provide now two dynamic programming principles. When the terminal reward $g$ is Borelian (resp. continuous), we show a {\em weak} (resp. {\em strong}) dynamic programming principle.
We stress that in the continuous case, the {\em strong} dynamic principle cannot be directly deduced from the weak one.

\subsection{A {\em weak} dynamic programming principle in the irregular case}

\paragraph{Measurability  properties of the functions $u$ and $u^{ \alpha }$.}
\noindent We first show some measurability properties of the functions $u^{ \alpha }(t,x)$ with respect  to 
control $\alpha$ and initial condition $x$.

By the a priori estimates on doubly reflected BSDEs provided in $\cite{DQS2}$ (see Proposition 5.3), we derive the following measurability property of $u^\alpha.$
\begin{lemma}[A measurability property of $u^{ \alpha }$] \label{mesu}
Let $s \in [0,T]$.

 The map $(\alpha,x) \mapsto u^{ \alpha } (s,x)$; $({\mathcal{A}}_s^s \times \mathbb{R}\,, \, 
 \mathcal{B}'({\mathcal{A}}_s^s)
 \otimes
 \mathcal{B}(  \mathbb{R}) ) \rightarrow ( \mathbb{R},  \mathcal{B}(  \mathbb{R}) ) $ is measurable.  
Here $ \mathcal{B}'({\mathcal{A}}_s^s)$ denotes the $\sigma$-algebra 
 induced by  ${\cal B}({\mathbb H}_s^2)$ on ${\mathcal{A}}_s^s$.  
\end{lemma}

\dproof 
Let $x^1, x^2$ $\in \mathbb{R}$, 
and $\alpha^1,\alpha^2 \in \mathcal{A}_s ^s.$
By classical results on diffusion processes,  we have
\begin{equation}\label{xesti}
 \mathbb{E}[\sup_{r \geq s }|X_r^{\alpha^1,s,x^1}-X_r^{\alpha^2,s,x^2}|^2]\leq 
 C (\| \alpha^1-\alpha^2 \|^2_{\mathbb{H}_s^2} + |x^1- x^2|^2 ). 
 \end{equation}
Let  ${Y}_{s,T}^{\alpha,s,x} [ \xi_., \zeta_., \eta]$ be the solution at time $s$
  of the  doubly reflected BSDE associated with driver $f^{\alpha,s,x}:= (f(\alpha_r,r, X_{r}^{\alpha,s,x},.){\bf 1}
  _{r \geq s}
  ) $, barriers 
  $\xi_r \leq \zeta_r$, $s\leq r <T$, terminal condition $\eta$. We here suppose that 
  ${\bar \xi}_r:= \xi_r {\bf 1}_{ r <T} + \eta {\bf 1}_{ r =T}$ and ${\bar \zeta}_r:= \zeta_r {\bf 1}_{ r <T} + \eta {\bf 1}_{ r =T}$ satisfy Mokobodzki's condition (see \eqref{mo} or Definition 3.9 in \cite{DQS2}). \\
  Using the Lipschitz property of 
$f$ with respect to $x, \alpha$, we have 
$$\| \sup_{y,z}|f(\alpha^1,X^{\alpha^1,s,x^1},y,z)-f(\alpha^2,X^{\alpha^2,s,x^2},y,z)| \,\|_{\mathbb{H}^2} \leq C(\|\alpha^1-\alpha^2\|_{\mathbb{H}^2}+\|X^{\alpha^1,s,x^1}-X^{\alpha^2,s,x^2}\|_{\mathcal{S}^2}).$$
  By the estimates on doubly reflected BSDEs with universal constants provided in \cite{DQS2} 
  (see Proposition 5.3), 
 we obtain that  
for all $x^1, x^2$ $\in \mathbb{R}$, 
$\alpha^1,\alpha^2 \in \mathcal{A}_s ^s$, barriers $\xi_.^1,\zeta_.^1\in S_s ^2$, 
$\xi_.^2,\zeta_.^2 \in S_s ^2$, and terminal conditions $\eta^1$, $\eta^2 \in L_s^2$,
\begin{align}\label{phi}
& |Y_s^{\alpha^1,s,x^1}[\xi^1_., \zeta_.^1, \eta^1]-Y_s^{\alpha^2,s,x^2}[\xi^2_., \zeta_.^2 , \eta^2]|^2  \nonumber\\
 &\leq 
 C' (\| \alpha^1-\alpha^2 \|^2_{\mathbb{H}_s^2} + |x^1- x^2|^2 +  \| \xi_.^1-\xi_.^2 \|^2_{S_s^2}+ \| \zeta_.^1-\zeta_.^2 \|^2_{S_s^2}+ \| \eta^1 - \eta^2\|^2_{L_s^2} ), 
 \end{align}
 where 
 $C'$ is a constant depending only on $T$ and the Lipschitz constant $C$ of the map $f$.
 
 Let $\Phi: \mathcal{A}_s ^s \times \mathbb{R} \times 
  ({\cal S}_s^2)^2 \times L_s^2 \rightarrow {\cal S}_s^2$; 
 $(\alpha,x, \xi_. , \zeta_., \eta) \mapsto {Y}_{s,T}^{\alpha,s,x} [ \xi_., \zeta_., \eta]$.
By \eqref{phi}, the map $\Phi$  is Lipschitz-continuous with respect to the norm 
  $\| \,. \,\|^2_{\mathbb{H}_s^2} + |\,.\,|^2  + \| \,. \,\|^2_{{\cal S}_s^2} + \| \,. \,\|^2_{{\cal S}_s^2}$ $+   \| \,. \,\|^2_{L_s^2}$.\\ 
 Now, by Proposition 3.6 in \cite{DQS3}, the map $L_s^{2p} \cap L_s^2 \rightarrow L_s^2$,  $\xi \mapsto g(\xi)$ is Borelian.
Hence, the map 
  $(\alpha,x) \mapsto (\alpha, x,h_1(.,X_{.}^{\alpha,s,x}),  h_2(.,X_{.}^{\alpha,s,x}) , g(X_{T}^{\alpha,s,x} ) )$   
  defined on $\mathcal{A}_s ^s 
   \times \mathbb{R} $ and 
  valued in 
  $\mathcal{A}_s ^s 
   \times \mathbb{R}  \times ({\cal S}_s^2)^2 \times L^2_s$ is 
${\cal B}(\mathcal{A}_s ^s ) \otimes {\cal B}(\mathbb{R}) /$ ${\cal B}(\mathcal{A}_s ^s )    \otimes {\cal B}(\mathbb{R}) \otimes {\cal B}(({\cal S}_s^2)^2) \otimes {\cal B}(L_s^2) $-measurable. \\
  By composition, it follows that the map 
$(\alpha,x) \mapsto  {Y}_{s,T}^{\alpha,s,x}[g(X^{\alpha,t,x}_T)] = u^\alpha(s,x)$ is 
Borelian. 
\fproof

Using the a priori estimates on doubly reflected BSDEs provided in $\cite{DQS2}$ and standard arguments, we derive the following result.
\begin{lemma}\label{poly}
For each $t \in [0,T]$, for each $\alpha \in \mathcal{A}_t ^t$, the map $u^\alpha$ has at most polynomial growth with respect to $x$. The property still holds for the value function $u$.
\end{lemma}

We  introduce the 
upper semicontinuous envelope  $u^*$ and the lower semicontinuous envelope  $u_*$ of the value function $u$, defined 
by
$$
u^*(t,x) := \limsup_{(t',x') \rightarrow (t,x)} u(t',x'); \quad 
 u_*(t,x) := \liminf_{(t',x') \rightarrow (t,x)} u(t',x')  \quad \forall (t, x) \in [0,T] \times {\mathbb R}.$$
 We also define the maps $\bar u^*$ and $\bar u_*$ for each $(t, x) \in [0,T] \times {\mathbb R}$ by
 $$
\bar u^*(t,x) :=  u^*(t,x) {\bf 1}_{t < T} + g(x) {\bf 1}_{t =  T}; \quad
\bar u_*(t,x) := u_*(t,x) {\bf 1}_{t < T} + g(x) {\bf 1}_{t =  T}. 
$$
 The functions $\bar u^*$ and $\bar u_*$ are Borelian. We have  $\bar u_* \leq u \leq \bar u^*$ and $\bar u_* (T,.)= u (T,.)=  \bar u^*(T,.)= g(.)$. Note that $ \bar u^*$ (resp. $\bar u_*$)  is not necessarily 
  upper (resp. lower) semicontinuous on $ [0,T] \times {\mathbb R}$, since the  terminal reward  $g$ is only Borelian.

 For each $ \theta \in {\cal T}$ and each  $\xi$ in $ L^2({\cal F}_ \theta)$, we denote by 
 $(Y^{\alpha,t,x}_{., \theta}(\xi), Z^{\alpha,t,x}_{., \theta}(\xi), k^{\alpha,t,x}_{., \theta}(\xi))$ 
 the unique solution in $\mathcal{S}^2 \times \mathbb{H}^2 \times
\mathbb{H}^2_{\nu}$ of the doubly reflected BSDE with  driver $f^{\alpha,t,x}{\bf 1}_ { \{ s \leq  \theta  \}  }$,
 terminal time $T$, 
 terminal condition $\xi$ and barriers 
 $h_i(r, X_r^{\alpha,t,x}) {\bf 1}_{r < \theta} + \xi {\bf 1}_{r \geq \theta}$.


 
 We now state a weak dynamic programming principle for our mixed game problem, which can be seen  as the analogous of the one  shown in \cite{DQS3} for a mixed optimal stopping/stochastic control problem.

\begin{theorem}\label{IMP}[Weak dynamic programming principle]
Suppose the subset ${\bf A}$ of ${\mathbb R}$ is compact.
We have the following {\em sub--optimality principle of dynamic programming}:\\
for each $t \in [0,T]$ and for each stopping time $\theta \in \mathcal{T}^t_{t},$
\begin{equation}\label{DPP}
u(t,x) \leq {\sup_{\alpha \in \mathcal{A}_t^t}}
Y_{t, \theta }^{\alpha,t,x}\left[{\bar u}^*(\theta,X_{\theta}^{\alpha,t,x})\right].
\end{equation}
We also have the  {\em super--optimality principle of dynamic programming}:\\
for each $t \in [0,T]$ and for each stopping time $\theta \in \mathcal{T}^t_{t},$
\begin{equation}\label{DPP2}
u(t,x) \geq {\sup_{\alpha \in \mathcal{A}_t^t}}Y_{t, \theta}^{\alpha,t,x}\left[{\bar u}_*(\theta,X_{\theta}^{\alpha,t,x})\right].
\end{equation}
These results still hold with $\theta$ replaced by $\theta^{\alpha}$ in inequalities \eqref{DPP} and \eqref{DPP2}, given a family of stopping times indexed by controls $\{ \theta^{\alpha}, \alpha \in \mathcal{A}_t^t\}.$

\end{theorem}

\dproof Without loss of generality, we can suppose that $t=0$.
Let us show  inequality \eqref{DPP}. For each $(t,x) \in$ $[0,T] \times \R$, we set
$\bar h_i(t, x):= h_i(t, x) {\bf 1}_{t < T} + g(x) {\bf 1}_{t =T}$.
Let $\theta \in \mathcal{T}$.
For each $n \in \mathbb N$, we define 
$
\theta^n:=\sum_{k=0}^{2^n-1}t_k\textbf{1}_{A_k}+T \textbf{1}_{\theta=T},
$
where  $t_k:=\frac{(k+1)T}{2^n}$ and $A_k:= \{\frac{kT}{2^n} \leq \theta <\frac{(k+1)T}{2^n}\}$. Note that 
$\theta^n \in \mathcal{T}$ and $\theta^n \downarrow \theta$.\\
Let $n\in \mathbb N$. Let $\alpha \in {\cal A}$.
By 
 the flow property for doubly reflected BSDEs, we get 
$
Y_{0,T}^{\alpha,0,x}= 
Y_{0,\theta^n}^{\alpha,0,x}[Y_{\theta^n,T}^{\alpha,\theta^n,X_{\theta^n}^{\alpha,0,x}}].
$
For each $s \in [0,T]$, for each $\omega \in \Omega$, set $^s\omega:= (\omega_{r \wedge s})_{0 \leq r \leq T}$. 
Note that for each $ \omega \in \Omega$, we have
$\omega =  \, ^s\omega +  \omega^s {\bf 1}_{]s,T]}$.
In the sequel, we identify $\omega$ with $(^s\omega, \omega^s).$

By the splitting property for doubly reflected BSDEs (see Proposition \ref{egal}), there exists a $P$-null set ${\cal N}$ such that for each $k$ and 
for each $\omega \in A_k \cap {\cal N}^c$, 
we have
$$Y_{\theta^n,T}^{\alpha,\theta^n,X_{\theta^n}^{\alpha,0,x}}(\omega)= 
Y_{t_k,T}^{\alpha,t_k,X_{t_k}^{\alpha,0,x}} (^{t_k}\omega)= 
Y_{t_k,T}^{\alpha
(^{t_k}\omega,\cdot),t_k,X_{t_k}^{\alpha,0,x}(^{t_k}\omega)} = u ^{\alpha (^{t_k}\omega,\cdot)}(t_k, X_{t_k}^{\alpha,0,x}(^{t_k}\omega)),$$
where the last equality follows from the definition of $u ^{\alpha (^{t_k}\omega,\cdot)}$.
Now, by definition of $u$, we have  $u ^{\alpha (^{t_k}\omega,\cdot)}(t_k, X_{t_k}^{\alpha,0,x}
(^{t_k}\omega)) \leq u (t_k, X_{t_k}^{\alpha,0,x}
(^{t_k}\omega))$. Since $u \leq \bar u^*$, we derive that
$Y_{\theta^n,T}^{\alpha,\theta^n,X_{\theta^n}^{\alpha,0,x}} \leq \bar u ^*(\theta^n,X_{\theta^n}^{\alpha,0,x})$ 
a.s.\, Hence, using the comparison theorem for doubly reflected BSDEs, we obtain
$ Y_{0,T}^{\alpha,0,x}= 
Y_{0,\theta^n}^{\alpha,0,x}[Y_{\theta^n,T}^{\alpha,\theta^n,X_{\theta^n}^{\alpha,0,x}}]  \leq  Y_{0,\theta^n}^{\alpha,0,x}[\bar u ^*(\theta^n,X_{\theta^n}^{\alpha,0,x})].$
By taking the limsup, we thus get 
\begin{equation}\label{limitesup}
 Y_{0,T}^{\alpha,0,x} \leq \lim \sup_{n \rightarrow \infty} Y_{0,\theta^n}^{\alpha,0,x}[\bar u ^*(\theta^n,X_{\theta^n}^{\alpha,0,x})].
 \end{equation}

Let us now show that the assumptions of the Fatou lemma for  
doubly reflected BSDEs (see Lemma \ref{fatou2} and Remark \ref{limsup}) are satisfied. 
Recall  that $u \leq {h_2}$ on $[0,T[ \times \R$. Since $h_2$ is continuous, we get $u^* \leq h_2$ on 
$[0,T[ \times \R$. 
On $\{ \theta<T \}$, we thus have 
 $$\lim \sup_{n \rightarrow \infty} \bar u ^*(\theta^n, X_{\theta^n}^{\alpha,0,x}) \leq \lim \sup_{n \rightarrow \infty} \overline{h}_2(\theta^n, X_{\theta^n}^{\alpha,0,x}) =h_2(\theta, X_{\theta}^{\alpha,0,x}) \quad \text{ a.s. }$$
where the last equality follows from the continuity property of $h_2$ on $[0,T[ \times \mathbb{R}$.\\
Now, on $\{ \theta=T \}$,   $\theta^n=T$. Hence, we have
$$\bar u ^*(\theta^n, X_{\theta^n}^{\alpha,0,x})= \bar u ^*(T, X_{T}^{\alpha,0,x})=g(X_{T}^{\alpha,0,x})=\bar{h_2}(T, X_{T}^{\alpha,0,x}).  $$
We thus get $\limsup_{n \rightarrow +\infty} \bar u ^*(\theta^n, X_{\theta^n}^{\alpha,0,x}) 
\leq \bar{h_2}(\theta, X_\theta^{\alpha,0,x}) \text{ a.s.}$ \\ Similarly, one can show that 
$\limsup_{n \rightarrow +\infty} \bar u ^*(\theta^n, X_{\theta^n}^{\alpha,0,x}) 
\geq \bar{h}_1(\theta, X_\theta^{\alpha,0,x}) \text{ a.s.}$ We thus have 
$$\bar{h}_1(\theta, X_\theta^{\alpha,0,x})  \leq \limsup_{n \rightarrow +\infty} \bar u _*(\theta^n, X_{\theta^n}^{\alpha,0,x}) 
\leq \bar{h}_2(\theta, X_\theta^{\alpha,0,x}) \text{ a.s.}$$
Condition \eqref{ij} (with limsup instead of liminf) is thus satisfied with $\xi^n = \bar u _*(\theta^n, X_{\theta^n}^{\alpha,0,x})$ and 
$\eta_t= \bar{h_1}(t, X_t^{\alpha,0,x})$ and $\zeta_t= \bar{h_2}(t, X_t^{\alpha,0,x})$. 

 Since $h_1$ is ${\cal C}^{1,2}$, by Proposition \ref{Moko}, 
 for each $n$, 
 $h_1(s,X_s^{\alpha,t,x}){\bf 1}_{\{s< \theta^n\}} + \xi^n {\bf 1}_{\{s= \theta^n\}}$ and \\ 
 $h_2(s,X_s^{\alpha,t,x}){\bf 1}_{\{s< \theta^n\}} + \xi^n {\bf 1}_{\{s= \theta^n\}}$ satisfy
  Mokobodzki 's condition. 
  
  Moreover, setting  $ \bar \xi := \limsup_{n \rightarrow +\infty} \bar u ^*(\theta^n, X_{\theta^n}^{\alpha,0,x}) $, the processes 
 $h_1(s,X_s^{\alpha,t,x}){\bf 1}_{\{s< \theta\}} + \bar \xi {\bf 1}_{\{s= \theta\}}$ and 
 $h_2(s,X_s^{\alpha,t,x}){\bf 1}_{\{s< \theta\}} + \bar \xi{\bf 1}_{\{s= \theta\}}$ also satisfy Mokobodzki 's condition.
 We can thus apply the  Fatou lemma for  
doubly reflected BSDEs (Lemma \ref{fatou2}).
Using \eqref{limitesup}, we then get:
$$ Y_{0,T}^{\alpha,0,x} \leq \lim \sup_{n \rightarrow \infty} Y_{0,\theta^n}^{\alpha,0,x}[\bar u ^*(\theta^n,X_{\theta^n}^{\alpha,0,x})]  \leq Y_{0,\theta}^{\alpha,0,x}[\lim \sup_{n \rightarrow \infty} \bar u ^*(\theta^n,X_{\theta^n}^{\alpha,0,x})].$$
Using the upper semicontinuity property of $\bar u ^*$ on $[0,T[ \times \R$ and  $\bar u ^*(T,x)= g(x)$, we  obtain
$$Y_{0,T}^{\alpha,0,x} \leq  Y_{0,\theta}^{\alpha,0,x}[\lim \sup_{n \rightarrow \infty} \bar u ^*(\theta^n,X_{\theta^n}^{\alpha,0,x})] \leq Y_{0,\theta}^{\alpha,0,x}[\bar u ^*(\theta,X_{\theta}^{\alpha,0,x})],$$

\noindent and this holds for each $\alpha \in \mathcal{A}$. Taking the supremum over 
$\alpha \in \mathcal{A}$,
 we get  \eqref{DPP}.\\
The proof of inequality \eqref{DPP2} relies on similar arguments as above as well as on an
existence result of $\varepsilon$-optimal controls satisfying appropriate measurability properties (see Proposition \ref{mesu2}). 
For more details on this last argument, we refer to the proof of the weak dynamic DPP for 
a mixed optimal stopping/control problem with ${\cal E}^f$-expectations 
(Th. 3.9 in \cite{DQS3}).
  \fproof

\subsection{ A {\em strong} dynamic programming principle in the continuous case} \label{continuous}
In this section, the set ${\bf A}$, where the controls are valued, is a nonempty closed 
subset of ${\mathbb R}^p$.

\noindent We suppose here that Assumption \ref{H2} holds. 


\begin{assumption}\label{H2} 
The terminal reward map $g$ is Lipschitz-continuous and the barriers maps $h_1$, $h_2$
 are Lipschitz-continuous with respect to $x$ uniformly in $t$. Moreover, we have 
 \begin{equation} \label{ca}
 h_1(T,x) \leq g(x) \leq h_2(T,x),\,\, \, \forall x \in {\mathbb R}.
 \end{equation}
 \end{assumption}
 
 Under this assumption, we show that $u$ is continuous with respect to $x$, and that $u^{\alpha}$ is continuous  with respect to $(\alpha,x)$.
\begin{lemma}[A continuity property of $u^{ \alpha }$ and $u$]\label{mesura} 

Suppose Assumption \ref{H2} holds. Then for each $s \in [0,T]$, the map  $(\alpha,x) \mapsto u^\alpha(s,x)$ is continuous 
and the value function $x \mapsto u(s,x)$ is continuous.
\end{lemma}

\dproof 
By estimates \eqref{xesti} and \eqref{phi}, we derive that 
\begin{align*}
|u^{\alpha^1} (s,x^1)-u^{\alpha^2} (s,x^2) |^2 &=     |Y_s^{\alpha^1,s,x^1}[g(X^{\alpha^1,t,x^1}_T)]-Y_s^{\alpha^2,s,x^2}[g(X^{\alpha^2,t,x^2}_T)]|^2 \\ & \leq 
 C' (\| \alpha^1-\alpha^2 \|^2_{\mathbb{H}_s^2} + |x^1- x^2|^2 ), 
 \end{align*}
 where 
 $C'$ is a constant depending only on $T$ and the Lipschitz constant $C$ of the map $f$.
 It follows that the map  $(\alpha,x) \mapsto u^\alpha(s,x)$ is Lipschitz-continuous.\\
  Hence, since $u(s,x)= \sup_{\alpha} u^{\alpha}(s,x)$, the map
 $x \mapsto u(s,x)$ is Lipschitz-continuous.
\fproof



In order to show the {\em strong} dynamic programming principle, we first prove that the value function $u$ is 
continuous with respect to $(t,x)$. We have already shown that $u$ is continuous with respect to
 $x$ uniformly in $t$ (see Lemma 
\ref{mesura}). 
%
It is thus sufficient to show the continuity of $u$ with respect to $t$. To this purpose, we first prove a {\em strong} dynamic programming principle at deterministic times.
\begin{lemma}\label{IMP1} Suppose that $g$, $h_1$ and $h_2$ satisfy the continuity Assumption \ref{H2}. Let $t \in [0,T]$.\\
For  all $s \geq t$, the value function $u$ defined by 
$\eqref{probform1}$ satisfies the equality 
\begin{equation}\label{DPPP}
u(t,x) = \sup_{\alpha \in \mathcal{A}_t^t}Y_{t,s}^{\alpha,t,x}\left[u(s,X_{s}^{\alpha,t,x})\right].
\end{equation}
\end{lemma}
\dproof 
We first show that:
\begin{equation}\label{sense1}
u(t,x) \leq \sup_{\alpha \in \mathcal{A}_t^t}Y_{t,s}^{\alpha,t,x}\left[u(s,X_{s}^{\alpha,t,x})\right].
\end{equation}
By the flow property for doubly reflected BSDEs (see \cite{Buckdahn2}), we have that:
$$
Y_{t,T}^{\alpha,t,x}= Y_{t,s}^{\alpha,t,x}[Y_{s,T}^{\alpha,t,x}].
$$
Note  that for almost-every $\omega$, at fixed $^s \omega$,  
the process $\alpha(^{s}\omega, T^s)$ (denoted also by $\alpha(^{s}\omega, \cdot)$ belongs to 
$\mathcal{A}_{s}^{s}$ (see Section \ref{A1} in the Appendix).
Moreover, by Proposition \ref{egal}, for almost-every $\omega$, at fixed $^{s}\omega$, $Y_{s,T}^{\alpha,t,x}(^{s}\omega)$ coincides with $Y_{s,T}^{\alpha(^{s}\omega, \cdot),s,X_{s}^{\alpha,t,x}(^{s}\omega)}$, the solution at time $s$ of the doubly reflected BSDE associated with control $\alpha(^{s}\omega, \cdot)$, initial conditions $s,X_{s}^{\alpha,t,x}(^{s}\omega)$, with respect to the filtration ${\mathbb F}^s$ and driven by the $s$-translated Brownian motion and $s$-translated Poisson measure.
Now,  by using the definition  of $u^\alpha$, we get that for almost-every $\omega$,
$$Y_{s,T}^{\alpha,t,x}(^{s}\omega)=Y_{s,T}^{\alpha(^{s}\omega, \cdot),s,X_{s}^{\alpha,t,x}(^{s}\omega)}=u^{\alpha({^s}\omega, \cdot)}(s,X_{s}^{\alpha,t,x}({^s}\omega)) \leq u(s,X_{s}^{\alpha,t,x}({^s}\omega)).$$ 
%
Finally, the comparison theorem for doubly reflected BSDEs (see  Th. 5.1 in \cite{DQS2})  leads to:
$$Y_{t,T}^{\alpha,t,x}= Y_{t,s}^{\alpha,t,x}[Y_{s,T}^{\alpha,t,x}] \leq Y_{t, s}^{\alpha,t,x}
[u(s,X_{s}^{\alpha,t,x})].$$
Taking the supremum over $\alpha \in \mathcal{A}_t^t$ in this inequality, we get inequality $\eqref{sense1}$.

It remains to show the   inequality:
\begin{equation} \label{equu1}
\sup_{\alpha \in \mathcal{A}_t^t}Y_{t,s}^{\alpha,t,x}\left[u(s,X_{s}^{\alpha,t,x})\right] \leq u(t,x).
\end{equation}


\noindent Fix $s \in [t,T]$ and $\alpha \in \mathcal{A}_t^t.$
 By  Corollary \ref{selectbis}, there exists an ``optimizing" sequence of measurable controls for $u(s,X_{s}^{\alpha,t,x})$, satisfying appropriate measurable properties. More precisely, 
there exists a  sequence $(\alpha^n)_{n \in {\mathbb N}}$ of controls belonging to 
$\mathcal{A}_{s}^{t}$ such that, for $P$-almost every $\omega$, we have
 \begin{equation}\label{310}
u(s,X_{s}^{\alpha,t,x}(^s \omega))=\lim_{n \rightarrow \infty } u^{\alpha^n(^s \omega, \cdot)}
( s, X_{s}^{\alpha,t,x}(^s \omega))= \lim_{n \rightarrow \infty } Y_{s,T}^{\alpha^{n}(^s \omega, \cdot),
s, X_{s}^{\alpha,t,x}(^s \omega)}.
\end{equation}
For each $n \in {\mathbb N}$, we set:
\begin{equation*}\label{equ2}
\Tilde{\alpha}^n_u:=\alpha_u \textbf{1}_{u<s}+\alpha^n_u\textbf{1}_{s \leq u \leq T}.
\end{equation*}
\noindent Note that $\Tilde{\alpha}^n \in \mathcal{A}_t^t$.
By the splitting property for doubly reflected BSDEs (see Proposition \ref{egal}), 
for $P$-almost every $\omega$, we have 
 \begin{equation*}\label{equ1}
Y_{s,T}^{\alpha^{n}(^s \omega, \cdot),
s, X_{s}^{\alpha,t,x}(^s \omega)}= Y_{s,T}^{\tilde \alpha^{n},t, x}(\omega).
\end{equation*}
Hence, by \eqref{310}, applying the continuity property of doubly reflected BSDEs with respect to terminal condition (see 
 Proposition \ref{continuity}), we obtain
 \begin{equation}\label{equ5}
Y_{t,s}^{\alpha,t,x}\left[u(s,X_{s}^{\alpha,t,x})\right] 
=Y_{t,s}^{\alpha,t,x}\left[\lim_{n \rightarrow \infty } Y_{s,T}^{\tilde \alpha^{n},t,x}\right] 
= \lim_{n \rightarrow \infty} Y_{t,s}^{\alpha,t,x}\left[ Y_{s,T}^{\tilde \alpha^{n},
t,x}\right].
\end{equation}
Now, by the flow property of doubly reflected BSDEs, for each $n$, we have $Y_{t,s}^{\alpha,t,x}\left[ Y_{s,T}^{\tilde \alpha^{n},
t,x}\right]= Y_{t,T}^{\Tilde{\alpha}^n,t,x}$.
We thus  get 
$$ Y_{t,s}^{\alpha,t,x}\left[u(s,X_{s}^{\alpha,t,x})\right] = \lim_{n \rightarrow \infty} Y_{t,T}^{\Tilde{\alpha}^n,t,x}\leq u(t,x).$$
Taking the supremum on $\alpha \in \mathcal{A}_t^t$ in this inequality, we get $\eqref{equu1}$, 
which ends the proof. 
\fproof

Using this {\em strong} dynamic programming principle at deterministic times, we derive the continuity property of $u$ with respect to time $t$.
\begin{theorem}\label{continu} Suppose that $g$, $h_1$ and $h_2$ satisfy the continuity Assumption \ref{H2}.
The value function  $u$ is then continuous with respect to $t$, uniformly in $x$.
\end{theorem}

\dproof
Since  \eqref{ca} holds, we have $h_1(t,x)\leq  u(t,x)\leq h_2(t,x)$ for all 
$(t,x) \in [0,T] \times {\mathbb R}$. \\
Let $0 \leq t < s \leq T$.
We have  
\begin{equation*} 
|u(t,x)-u(s,x)| \leq |u(t,x)-\sup_{\alpha \in \mathcal{A}_t^t}\mathcal{E}_{t,s}^{\alpha,t,x}[u(s, X_{s}^{\alpha,t,x})]|+|\sup_{\alpha \in \mathcal{A}_t^t}\mathcal{E}_{t,s}^{\alpha,t,x}[u(s, X_{s}^{\alpha,t,x})]-u(s,x)|.
\end{equation*}
We start by estimating $|\sup_{\alpha \in \mathcal{A}_t^t}\mathcal{E}_{t,s}^{\alpha,t,x}[u(s, X_{s}^{\alpha,t,x})]-u(s,x)|.$
\begin{align}
&| \sup_{\alpha \in \mathcal{A}_t^t}\mathcal{E}_{t,s}^{\alpha,t,x}[u(s, X_{s}^{\alpha,t,x})]-u(s,x)| \leq \sup_{\alpha \in \mathcal{A}_t^t}|\mathcal{E}_{t,s}^{\alpha,t,x}[u(s, X_{s}^{\alpha,t,x})]-u(s,x)| \nonumber \\
&\leq \sup_{\alpha \in \mathcal{A}_t^t}|\mathcal{E}_{t,s}^{\alpha,t,x}[u(s, X_{s}^{\alpha,t,x})]-\mathcal{E}_{t,s}^{0}[u(s,x)]| \leq C {\mathbb E}[\sup_{t \leq r \leq s}(X_{r}^{\alpha,t,x}-x)^2]^{\frac{1}{2}}  \leq C|s-t|(1+x^{2p})^{\frac{1}{2}}.
\label{ineq1}
 \end{align}
Here, $\mathcal{E}^0$ denotes the conditional expectation associated with the driver equal to $0$.
In order to obtain the above relation, we have used BSDEs estimates (see \cite{16}), the Lipschitz property
 of $u$ with respect to $x$ (see Lemma \ref{mesu}) and the polynomial growth of $u$ (see Lemma \ref{poly}). \\
We now estimate $|u(t,x)-\sup_{\alpha \in \mathcal{A}_t^t}\mathcal{E}_{t,s}^{\alpha,t,x}[u(s, X_{s}^{\alpha,t,x})]|.$
Using the strong dynamic programming principle for deterministic times (see Lemma \ref{IMP1}), we derive that:
$$u(t,x) = \sup_{\alpha \in \mathcal{A}_t^t}Y_{t,s}^{\alpha,t,x}\left[u(s,X_{s}^{\alpha,t,x})\right].$$
Now, the solution $Y$ of the doubly reflected BSDE is smaller than the solution of the reflected BSDE with the same lower barrier, denoted by $\Tilde{Y}$. Hence,
$$Y_{t,s}^{\alpha,t,x}\left[u(s,X_{s}^{\alpha,t,x})\right] \leq \Tilde{Y}_{t,s}^{\alpha,t,x}\left[h_1(r,X_{r}^{\alpha,t,x})\textbf{1}_{t \leq r < s}+u(s,X_{s}^{\alpha,t,x})\textbf{1}_{r=s}\right].$$
By the characterization of the solution of a reflected BSDE, we derive that
\begin{eqnarray}\label{inegalite1}
 u(t,x)-  \sup_{\alpha \in \mathcal{A}_t^t}\mathcal{E}_{t,s}^{\alpha,t,x}[u(s, X_{s}^{\alpha,t,x})] 
 &\leq& \sup_{\alpha \in \mathcal{A}_t^t} \sup_{\tau \in \mathcal{T}^t_{t}}(\mathcal{E}_{t, \tau \wedge s}^{\alpha,t,x}[h_1(\tau,X_{\tau}^{\alpha,t,x})\textbf{1}_{\tau<s}+u(s, X_{s}^{\alpha,t,x})\textbf{1}_{\tau \geq s}]
 \nonumber  \\
 &&  -\mathcal{E}_{t,s}^{\alpha,t,x}[u(s, X_{s}^{\alpha,t,x})\textbf{1}_{\tau < s}+u(s, X_{s}^{\alpha,t,x})\textbf{1}_{\tau \geq s}])\nonumber\\
& \leq & A, 
\end{eqnarray}
where $$A:= \sup_{\alpha \in \mathcal{A}_t^t} \sup_{\tau \in \mathcal{T}^t_{t}}|\mathcal{E}_{t, s  }^{\Tilde{f}^\alpha}[h_1(\tau,X_{\tau}^{\alpha,t,x})\textbf{1}_{\tau<s}+u(s, X_{s}^{\alpha,t,x})\textbf{1}_{\tau \geq s}]-\mathcal{E}_{t, s}^{f^\alpha}[h_1(s, X_{s}^{\alpha,t,x})\textbf{1}_{\tau < s}+u(s, X_{s}^{\alpha,t,x})\textbf{1}_{\tau \geq s}]| $$
with $\Tilde{f}^\alpha(s, \cdot):=f^{\alpha,t,x}(s, \cdot)\textbf{1}_{s \leq \tau},$ because $u \geq h_1.$
By the Lipschitz property in $x$ of $h_1$, the polynomial growth of $h_1$ and $f$ in $x$, and the standard estimates for BSDEs and SDEs, we have
 \begin{eqnarray}\label{inequno}
A^2
& \leq &C \sup_{\alpha \in \mathcal{A}_t^t} \sup_{\tau \in \mathcal{T}^t_{t}}\mathbb{E}[(h_1(\tau,X_{\tau}^{\alpha,t,x})\textbf{1}_{\tau<s}-h_1(s, X_{s}^{\alpha,t,x})\textbf{1}_{\tau<s})^2]+\mathbb{E}[\int_{\tau \wedge s}^{s}f^2(\alpha_r,r,h_1(\tau,X_{\tau}^{\alpha,t,x}),0,0)dr]  \nonumber  \\
& \leq & C[\sup_{t \leq r < s} \sup_x |h_1(r,x)-h_1(s,x)|^2+(1+|x|^p)^2|s-t|+|s-t|(1+|x|^{q})]. 
 \end{eqnarray}
Also, the solution $Y^{\alpha,t,x}$ of the doubly reflected BSDE is greater than the solution of the reflected BSDE with the same upper barrier, denoted by $\bar{Y}^{\alpha,t,x}$. Hence,
$$Y_{t,s}^{\alpha,t,x}\left[u(s,X_{s}^{\alpha,t,x})\right] \geq \bar{Y}_{t,s}^{\alpha,t,x}\left[h_2(r,X_{r}^{\alpha,t,x})\textbf{1}_{t \leq r < s}+u(s,X_{s}^{\alpha,t,x})\textbf{1}_{r=s}\right].$$
By the characterization of the solution of a reflected BSDE with upper barrier , we derive that
$$\bar{Y}_{t,s}^{\alpha,t,x}\left[h_2(r,X_{r}^{\alpha,t,x})\textbf{1}_{t \leq r < s}+u(s,X_{s}^{\alpha,t,x})\textbf{1}_{r=s}\right]= \inf_{\sigma \in \mathcal{T}^t_{t}}\mathcal{E}_{t, \sigma \wedge s}^{\alpha,t,x}[h_2(\sigma,X_{\sigma}^{\alpha,t,x})\textbf{1}_{\sigma<s}+u(s, X_{s}^{\alpha,t,x})\textbf{1}_{\sigma \geq s}]$$
It follows that
\begin{eqnarray}\label{inegalite2}
 u(t,x)-  \sup_{\alpha \in \mathcal{A}_t^t}\mathcal{E}_{t,s}^{\alpha,t,x}[u(s, X_{s}^{\alpha,t,x})] 
 &\geq& \sup_{\alpha \in \mathcal{A}_t^t} \inf_{\sigma \in \mathcal{T}^t_{t}}(\mathcal{E}_{t, \sigma \wedge s}^{\alpha,t,x}[h_2(\sigma,X_{\sigma}^{\alpha,t,x})\textbf{1}_{\sigma<s}+u(s, X_{s}^{\alpha,t,x})\textbf{1}_{\sigma \geq s}] \nonumber\\
 &&  -\mathcal{E}_{t,s}^{\alpha,t,x}[u(s, X_{s}^{\alpha,t,x})\textbf{1}_{\sigma< s}+u(s, X_{s}^{\alpha,t,x})\textbf{1}_{\sigma \geq s}]) \nonumber\\
& \geq & B, 
\end{eqnarray}
where $$B:= \sup_{\alpha \in \mathcal{A}_t^t} \inf_{\sigma \in \mathcal{T}^t_{t}}|\mathcal{E}_{t, s  }^{\Bar{f}^\alpha}[h_2(\sigma,X_{\sigma}^{\alpha,t,x})\textbf{1}_{\sigma<s}+u(s, X_{s}^{\alpha,t,x})\textbf{1}_{\sigma \geq s}]-\mathcal{E}_{t, s}^{f^\alpha}[h_2(s, X_{s}^{\alpha,t,x})\textbf{1}_{\sigma < s}+u(s, X_{s}^{\alpha,t,x})\textbf{1}_{\sigma \geq s}]| $$
with $\Bar{f}^\alpha(s, \cdot):=f^{\alpha,t,x}(s, \cdot)\textbf{1}_{s \leq \sigma},$ because $u \leq h_2.$
Now, by the Lipschitz property in $x$ of $h_2$, the polynomial growth of $h_2$ and $f$ in $x$, we have
 \begin{eqnarray}\label{ineqdos}
B^2
& \leq &C \sup_{\alpha \in \mathcal{A}_t^t} \inf_{\sigma \in \mathcal{T}^t_{t}}\mathbb{E}[(h_2(\sigma,X_{\sigma}^{\alpha,t,x})\textbf{1}_{\sigma<s}-h_2(s, X_{s}^{\alpha,t,x})\textbf{1}_{\sigma<s})^2]+\mathbb{E}[\int_{\sigma \wedge s}^{s}f^2(\alpha_r,r,h_2(\sigma,X_{\sigma}^{\alpha,t,x}),0,0)dr]  \nonumber  \\
& \leq & C[\sup_{t \leq r < s} \sup_x |h_2(r,x)-h_2(s,x)|^2+(1+|x|^p)^2|s-t|+|s-t|(1+|x|^{q})]. 
 \end{eqnarray}
 By \eqref{inegalite1} and \eqref{inegalite2}, we get
$$
 \vert u(t,x)-  \sup_{\alpha \in \mathcal{A}_t^t}\mathcal{E}_{t,s}^{\alpha,t,x}[u(s, X_{s}^{\alpha,t,x})] \vert ^2 
 \leq A^2 + B^2.$$
This, together with inequalities \eqref{inequno}, \eqref{ineqdos} and inequality \eqref{ineq1} implies that
\begin{align*}
 |u(t,x)-u(s,x)|^2  \leq & \,C[\sup_{t \leq r < s} \sup_x (|h_1(r,x)-h_1(s,x)|^2+ \sup_{t \leq r < s} \sup_x (|h_2(r,x)-h_2(s,x)|^2 ]  \\
 & +C[(s-t)^2 +|s-t|((1+|x|^p)^2+ 1+|x|^{q}).
\end{align*}
Now, since $h_1$ and $h_2$  are continuous with respect to $t$ on $[0,T]$ uniformly in $x$, we derive that 
$u$ is  continuous with respect to time $t$, uniformly in $x$.
The proof is thus ended.
\fproof

From this theorem, we derive the continuity property of $u$ with respect to $(t,x)$.
\begin{corollary} \label{certain} Suppose that $g$, $h_1$ and $h_2$ satisfy the continuity Assumption \ref{H2}.\\
The value function $u$ is then continuous on $[0,T] \times {\mathbb R}$. 
\end{corollary}
\dproof
 Since Assumption \ref{H2} is satisfied, $u$ is continuous with respect to $x$ (see Lemma \ref{mesu}). This property together with the previous theorem implies that $u$ is continuous with respect to $(t,x)$.
 \fproof
 
This result yields that under Assumption \ref{H2}, we have $u^*=u_*=u$.
Hence, by the {\em weak} dynamic programming principle (which still holds even if the set ${\bf A}$ is not compact because of the continuity Assumption \ref{H2}) it follows that the value function $u$ satisfies the following  {\em strong} dynamic programming principle at stopping times. 
\begin{theorem}\label{IMP2}[Strong dynamic programming principle] Suppose that $g$, $h_1$ and $h_2$ satisfy the continuity Assumption \ref{H2}.
 For each $t \in [0,T]$ and for each stopping time $\theta \in \mathcal{T}^t_{t},$ we have
\begin{equation}\label{eq}
u(t,x) ={\sup_{\alpha \in \mathcal{A}_t^t}}
Y_{t, \theta }^{\alpha,t,x}\left[u(\theta,X_{\theta}^{\alpha,t,x})\right].
\end{equation}
This result still holds with $\theta$ replaced by $\theta^{\alpha}$, given a family of stopping times indexed by controls $\{ \theta^{\alpha}, \alpha \in \mathcal{A}_t^t\}.$

\end{theorem}

\section{Generalized HJB variational inequalities }\label{sec5}
In this section, we do not suppose that 
Assumption \ref{H2} holds.

\noindent We introduce the  following  Hamilton Jacobi Bellman variational inequalities  (HJBVIs):

\begin{equation}\label{edp}
\begin{cases}
\,\,h_1(t,x) \leq u(t,x) \leq h_2(t,x),\quad t<T\\
\text{ if } u(t,x)<h_2(t,x) \text{ then } \inf_{\alpha \in A} \mathcal{H}^{\alpha} u(t,x) \geq 0\\
\text{ if } h_1(t,x)<u(t,x) \text{ then } \inf_{\alpha \in A} \mathcal{H}^{\alpha} u(t,x) \leq 0,
\end{cases}
\end{equation}
with terminal condition $u(T,x)=g(x), x\in \mathbb{R}$.
Here, $L ^{\alpha}:=A^{\alpha}+K^{\alpha},$ and for $\phi \in C^2(\mathbb{R})$,
\begin{itemize}

\item
$A^{\alpha}\phi(x) := \dfrac{1}{2}\sigma^2(x,\alpha)\dfrac{\p^2 \phi}{\p x^2}(x)+ b(x,\alpha) \dfrac{\p \phi }{\p x}(x)$ 
\item
$K^{\alpha} \phi(x) :=\int_{{\bf E}}\left(\phi(x+ \beta(x,\alpha,e))-\phi(x)- \dfrac{\p \phi}{\p x}(x)\beta(x,\alpha,e)\right) \nu (de)$
\item
$B^{\alpha} \phi(x) :=\phi(x+\beta(x,\alpha,\cdot))-\phi(x).$
\item
$\mathcal{H}^{\alpha} \phi(x) :=-\dfrac{\p u}{\p t}(t,x)-L^{\alpha}u(t,x)-f(\alpha,t,x,u(t,x), (\sigma \dfrac{\p u }{\p x})(t,x),B^{\alpha}u(t,x)) .$
\end{itemize}
\begin{definition}\rm
$\bullet$ An upper semicontinuous  function $u$ is said to be a {\em viscosity subsolution} of \eqref{edp}  if for any point $(t_0,x_0) \in [0,T[ \times \mathbb{R}$ and for any $\phi \in C^{1,2}([0,T] \times \mathbb{R})$ such that $\phi(t_0,x_0)=u(t_0,x_0)$ and $\phi-u$ attains its minimum at $(t_0,x_0)$, if $h_1(t_0,x_0)<u(t_0,x_0)$ then $\inf_{\alpha \in A} \mathcal{H}^{\alpha} \phi(t,x) \leq 0$.

$\bullet$
 A lower semicontinuous function $u$  is said to be a {\em  viscosity supersolution } of \eqref{edp}  if for any point $(t_0,x_0) \in [0,T[ \times \mathbb{R}$ and for any $\phi \in C^{1,2}([0,T] \times \mathbb{R})$ such that $\phi(t_0,x_0)=u(t_0,x_0)$ and $\phi-u$ attains its maximum at $(t_0,x_0)$, if $u(t_0,x_0)<h_2(t_0,x_0)$ then $\inf_{\alpha \in A}  \mathcal{H}^{\alpha} \phi(t,x) \geq 0$.

%
%
\end{definition}


 
\subsection{The irregular case}

 Using the weak dynamic programming principle (Theorem \ref{IMP}), we now prove that the value function of our problem is a  viscosity solution of the above HJBVI in a weak sense.

\begin{theorem}\label{existprob} Suppose that ${\bf A}$ is compact. The map $u$ is  a {\em weak viscosity solution} of \eqref{edp} in the sense that $u^*$  is a viscosity subsolution of  \eqref{edp} and 
$u_*$  is a viscosity supersolution of  \eqref{edp}.
\end{theorem}

\begin{remark}
Using this theorem, when $h_1(T,x) \leq g(x) \leq h_2(T,x)$ for all $x \in \mathbb{R}$, we show in the next section that 
$u_*$  is a viscosity supersolution of  \eqref{edp} with terminal value greater than $g_*$ (see Corollary \ref{1}).
Moreover, when $g$ is l.s.c, we show that the value function $u$ 
is the minimal viscosity supersolution of  \eqref{edp} with terminal value greater than $g$ 
(see Theorem \ref{5}).
\end{remark}

\dproof
We first prove that $u_*$  is a viscosity supersolution of \eqref{edp}.\\
Let $(t_0,x_0) \in [0,T[ \times \mathbb{R}$ and $\phi \in C^{1,2}([0,T] \times \mathbb{R})$ be such that $\phi(t_0, x_0)=u_*(t_0,x_0)$  and \\
$\phi(t,x) \leq u_*(t,x)$, for each $ (t,x) \in [0,T] \times \mathbb{R}$. 
Without loss of generality, we can suppose that the maximum is strict in $(t_0,x_0)$.\\
Let us assume that $u_*(t_0,x_0)<h_2(t_0,x_0)$ and that
\begin{equation*}
\inf_{\alpha \in A}\left(-\dfrac{\p}{\p t}\phi(t_0,x_0)-L^\alpha\phi(t_0,x_0)-f \left(\alpha,t_0,x_0, \phi (t_0,x_0), (\sigma \dfrac{\p \phi}{\p x})(t_0,x_0),B^\alpha \phi(t_0,x_0)\right)\right) < 0.
\end{equation*}
By continuity, we can suppose that there exists $\alpha \in A$, $\varepsilon>0$ and $\eta_\varepsilon>0$ such that:\\
$\forall (t,x)$ such that $t_0 \leq t \leq t_0+ \eta_\varepsilon<T$ and $|x-x_0| \leq \eta_\varepsilon$, we have 
$\phi(t,x) \leq h_2(t,x) -\varepsilon$ and
\begin{equation}\label{boule1}
 -\dfrac{\p}{\p t} \phi(t,x)-L^\alpha \phi(t,x)-f\left(\alpha,t,x,\phi(t,x), (\sigma \dfrac{\p \phi }{\p x})(t,x),B^\alpha \phi(t,x)\right) \leq -\varepsilon.
\end{equation}
Let $B_{\eta_\varepsilon}(t_0,x_0)$  be the ball of radius $\eta_\varepsilon$ and center $(t_0,x_0)$. Let $(t_n,x_n)_{n \in {\mathbb N}}$ be a sequence in $B_{\eta_\varepsilon}(t_0,x_0)$  with $t_n \geq t_0$ for each $n$, such that the sequence $(t_n,x_n,u(t_n,x_n))_{n \in {\mathbb N}}$ tends to $ (t_0,x_0,u_*(t_0,x_0))$.
We introduce the state process $X^{\alpha,t_n,x_n}$ associated with the above constant control $\alpha$.\\
Let $\theta^{\alpha,n}$ be the stopping time defined by
\begin{equation*}\theta^{\alpha,n}:= (t_0 + \eta_\varepsilon) \wedge \inf\{s \geq t_n \, , \, |X_s^{\alpha,t_n,x_n}-x_0| \geq \eta_\varepsilon \}.
\end{equation*}
%
By applying It\^o's lemma to $\phi(s, X_s^{\alpha,t_n,x_n})$, one can derive 
  that  
  \begin{equation}\label{bi}
  \left(\phi(s,X_s^{\alpha,t_n,x_n}), (\sigma \dfrac{\p \phi}{\p x})(s, X_s^{\alpha,t_n,x_n}), B^\alpha\phi(s,X_{s^-}^{\alpha,t_n,x_n}); s \in [t_n, \theta^{\alpha,n}]\right)
  \end{equation}
   is the solution of the BSDE associated with 
terminal time $\theta^{\alpha,n}$,  terminal value $\phi(\theta^{\alpha,n}, X_{\theta^{\alpha,n}}^{\alpha,t_n,x_n})$ and driver $-\psi^\alpha(s,X_s^{\alpha,t_n,x_n})$, where $\psi^{\alpha}(s,x):=\dfrac{\p}{\p s}\phi(s,x)+L^{\alpha}\phi(s,x)$.
By the definition of the stopping time $\theta^{\alpha,n}$ together with inequality \eqref{boule1}, we obtain:
\begin{equation}\label{eqq2}
-\psi^{\alpha}(s,X_s^{\alpha,t_n,x_n}) \leq f(\alpha,s,X_s^{\alpha,t_n,x_n},\phi (s,X_s^{\alpha,t_n,x_n}),  (\sigma \dfrac{\p \phi}{\p x} )(s,X_s^{\alpha,t_n,x_n}),B^\alpha \phi (s,X_s^{\alpha,t_n,x_n}))- \varepsilon,
\end{equation}
for $t_n \leq s \leq \theta^{\alpha,n}$ $ds \otimes dP$-a.s.\, 
 Now, since the maximum $(t_0,x_0)$ is strict, there exists $\gamma_\varepsilon$, which depends on $\eta_\varepsilon$, such that 
\begin{equation}\label{ut}
u_*(t,x) \geq \phi(t,x) +\gamma_\varepsilon \text{  on  }  [0,T] \times \mathbb{R} \setminus B_{\eta_\varepsilon}(t_0,x_0). 
\end{equation}
Note now that
\begin{equation*}\label{wifi1}
\phi(\theta^{\alpha,n}\wedge t, X_{\theta^{\alpha,n} \wedge t}^{\alpha,t_n,x_n})=\phi(t, X_{t}^{\alpha,t_n,x_n})\textbf{1}_{t<\theta^{\alpha,n}}+\phi(\theta^{\alpha,n}, X_{\theta^{\alpha,n}}^{\alpha,t_n,x_n})\textbf{1}_{t \geq\theta^{\alpha,n}}, \,\,\, t_n \leq t \leq T.
\end{equation*}
 Using  the inequalities $\phi(t,x) \leq h_2(t,x) -\varepsilon$ and \eqref{ut} together with the definition of $\theta^{\alpha,n}$, we derive that for each $t \in [t_n, \theta^{\alpha,n}]$,
$$
\phi( t, X_{ t}^{\alpha,t_n,x_n}) \leq (h_2(t, X_{t}^{\alpha,t_n,x_n})-\delta_\varepsilon)\textbf{1}_{t<\theta^{\alpha,n}}+(u_*(\theta^{\alpha,n}, X_{\theta^{\alpha,n}}^{\alpha,t_n,x_n})-\delta_\varepsilon)\textbf{1}_{t = \theta^{\alpha,n}} \quad {\rm a.s.}\,,
$$
where $\delta_\varepsilon:=\min(\varepsilon, \gamma_{\varepsilon})$.
Hence, by the inequality \eqref{eqq2} between the driver process
 $-\psi^{\alpha}(s,X_s^{\alpha,t_n,x_n})$ and the driver $f(\alpha,s,X_s^{\alpha,t_n,x_n}, .)$ computed along the solution \eqref{bi}, applying 
a comparison theorem between a BSDE and a reflected BSDE (see Proposition A.11 in \cite{DQS3}), we obtain
\begin{equation}\label{do}
\phi(t_n,x_n) \leq \Bar{Y}_{t_n, \theta^{\alpha,n}}^{\alpha,t_n,x_n}[h_2(t, X_{t}^{\alpha,t_n,x_n})\textbf{1}_{t<\theta^{\alpha,n}}+u_*(\theta^{\alpha,n}, X_{\theta^{\alpha,n}}^{\alpha,t_n,x_n})\textbf{1}_{t=\theta^{\alpha,n}}]-\delta_\varepsilon K ,
\end{equation}
where $\Bar{Y}$ is the solution of the reflected BSDE associated with upper barrier\\
$h_2(t, X_{t}^{\alpha,t_n,x_n})\textbf{1}_{t<\theta^{\alpha,n}}+u_*(\theta^{\alpha,n}, X_{\theta^{\alpha,n}}^{\alpha,t_n,x_n})\textbf{1}_{t=\theta^{\alpha,n}}$ and $K$ is a positive constant which only depends on $T$  and the Lipschitz constant of $f$. \\
Moreover, since the solution of the reflected BSDE with upper barrier is smaller than the solution of a doubly reflected BSDE with the same upper barrier, we have
\begin{eqnarray*}
\Bar{Y}_{t_0, \theta^{\alpha,n}}^{\alpha,t_n,x_n}[h_2(t, X_{t}^{\alpha,t_n,x_n})\textbf{1}_{t<\theta^{\alpha,n}}+u_*(\theta^{\alpha,n}, X_{\theta^{\alpha,n}}^{\alpha,t_n,x_n})\textbf{1}_{t=\theta^{\alpha,n}}] & \leq& Y^{\alpha,t_n,x_n} _{t_n,\theta^{\alpha,n}}[u_*(\theta^{\alpha,n}, X_{\theta^{\alpha,n}}^{\alpha,t_n,x_n})].
\end{eqnarray*}
Hence, using inequality \eqref{do}, we get
\begin{equation}\label{contract3}
\phi(t_n,x_n) 
 \leq Y_{t_n,\theta^{\alpha,n}}^{\alpha,t_n,x_n}[u_*(\theta^{\alpha,n}, X_{\theta^{\alpha,n}}^{\alpha,t_n,x_n})]-\delta_\varepsilon K. 
\end{equation}
Now, $\phi$ is continuous 
with $\phi(t_0,x_0)=u_*(t_0,x_0)$, and the sequence $(t_n,x_n,u(t_n,x_n)) $ converges to $(t_0,x_0,u_*(t_0,x_0))$ as $n$ tends to $+ \infty$.
We can thus assume that $n$ is sufficiently large so that
$
|\phi(t_n,x_n)-u(t_n,x_n)| \leq \dfrac{\delta_\varepsilon K}{2}
$. By inequality \eqref{contract3}, we thus get
\begin{equation*}\label{contract}
u(t_n,x_n) \leq Y_{t_n, \theta^{\alpha,n}}^{\alpha,t_n,x_n}[u_*(\theta^{\alpha,n}, X_{\theta^{\alpha,n}}^{\alpha,t_n,x_n})]-\dfrac{\delta_\varepsilon K}{2}. 
\end{equation*}
Now, by the super-optimality dynamic programming principle \eqref{DPP2}, since ${\bar u}_* \geq u_*$, we have 
$u(t_n,x_n) \geq Y_{t_n, \theta^{\alpha,n}}^{\alpha,t_n,x_n}[u_*(\theta^{\alpha,n}, X_{\theta^{n,\alpha}}^{\alpha,t_n,x_n})].$
We thus obtain a contradiction. Hence, $u_*$  is a viscosity supersolution of \eqref{edp}. 
It remains to  prove that $u^*$  is a viscosity subsolution of \eqref{edp}. The proof, which is based on quite similar arguments, is omitted.
\fproof

\begin{remark}
In the case of a Brownian framework and a continuous terminal reward, this property corresponds to a result shown in \cite{Buckdahn2} by using a penalization approach.
\end{remark}

\subsection{The continuous case}

\begin{theorem}\label{existprob2} Suppose that the set ${\bf A}$ is a closed subset of ${\mathbb R}^p$ and that 
the continuity Assumption \ref{H2} holds.
Then the value function $u$ is a viscosity solution in the classical sense, that is both a viscosity sub- and super-solution of the generalized HJBVIs  \eqref{edp}.
\end{theorem}
\dproof
Since Assumption \ref{H2} holds, by Corollary \ref{certain}, the value function $u$ is continuous with respect to $(t,x)$, which implies that $u^*=u_*=u$. Moreover, by Theorems \ref{IMP} and \ref{existprob} (which do not require the compactness of ${\bf A}$ in the continuous case), $u$ is a {\em weak} viscosity
 solution of  the  generalized HJBVIs \eqref{edp}. The result thus follows.
 \fproof

\paragraph{An uniqueness result.}
Suppose Assumption \ref{H2} holds and ${\bf A}$ is compact.
We assume moreover that 
${\bf E } = {\mathbb R}^*$, ${\cal K} = {\cal B} ({\mathbb R}^*)$ and $\int_{\bf E }(1 \wedge e^2) \nu (de) < + \infty$, and the following
\begin{assumption} \label{9.8}
\begin{itemize}

\item[1.] $f(\alpha,x,y,z,k) :=\overline{f}(\alpha,x,y,z,\int_{{\bf E}}k(e)\gamma(x,e)\nu(de))$
where \\ $\overline{f}: A \times [0,T] \times \mathbb{R}^4 \rightarrow \mathbb{R}$ is Borelian and satisfies:
\begin{itemize}

\item[(i)]
 $|\overline{f}(\alpha,t,x,0,0,0)| \leq C,$  for any $x\in \mathbb{R}$, $t\in [0,T]$, $\alpha \in A$.

\item[(ii)]
$ |\overline{f}(\alpha,t,x,y,z,k)- \overline{f}(\alpha',t,x',y',z',k')| \leq C(|\alpha-\alpha'|+ |y-y'|+|z-z'|+|k-k'|)$,  for any 
$x,x' \in \mathbb{R}$, $t\in [0,T]$, $y, y' \in \mathbb{R}$, $z,z' \in \mathbb{R}$, $k,k' \in \mathbb{R}$, $\alpha , \alpha ' \in A.$

\item[(iii)]
$k \rightarrow \overline{f}(\alpha,t,x,y,z,k)$ is non-decreasing,  for any $(\alpha,t,x,y,z,k) \in  A \times [0,T] \times \mathbb{R}^4.$
\end{itemize}

\item[2.] For each $R>0$, there exists a continuous function $m_R: \mathbb{R}_{+} \rightarrow \mathbb{R}_+$   such that
$m_R(0)=0$ and $$|\overline{f}(\alpha,t,x,v,p,q)-\overline{f}(\alpha,t,y,v,p,q)| \leq m_{R}(|x-y|(1+|p|)),$$ for any $  t \in [0,T], |x|, |y|,|v|\leq R, p,q \in \mathbb{R}, \alpha \in A.$

\item[3.]
$\gamma: \mathbb{R} \times {\bf E} \rightarrow \mathbb{R}$ is $ \mathcal{B}(\mathbb{R}) \otimes \mathcal{B}(\bf{E})$-measurable,  \\
$|  \gamma(x,e)-\gamma(x',e)| < C|x-x'| (1 \wedge e^2), x, x' \in \mathbb{R}, e \in \bf{E}$\\
  $|\gamma(x,e)| \leq C(1 \wedge |e|)$ and $\gamma(x,e)\geq 0$, $e \in \bf{E}$

\item[4.]
\text{  There exists }$r>0$ such that for any $x\in \mathbb{R}$, $t\in [0,T]$, $u,v\in \mathbb{R}$, $p\in \mathbb{R}$, $l\in \mathbb{R}$, $\alpha \in A$:
$$
\overline{f}(\alpha,t,x,v,p,l) - \overline{f}(\alpha,t,x,u,p,l)\geq  r(u-v)
\text{ when } u\geq v.$$


\item[5.] $|g(x)| + |h(t,x)| \leq C,$ for any $x\in \mathbb{R}$, $t\in [0,T]$.
\end{itemize}
\end{assumption}

\begin{lemma}[Comparison principle]\label{8.9} Suppose that the above assumptions hold. If $U$ is a bounded viscosity subsolution  and $V$ is a bounded viscosity supersolution of the HJBVIs  \eqref{edp} with
$U(T,x) \leq g(x) \leq V(T,x)$, $x \in \mathbb{R}.$ Then, $U(t,x) \leq V(t,x)$,
for each $(t,x) \in [0,T] \times \mathbb{R}$. 
\end{lemma}
The proof is similar to the case studied in \cite{DQS2} without controls  and is thus omitted. 

Moreover, using this comparison principle together with Theorem \ref{existprob2}, we derive the
 characterization of $u$ as the unique solution of the HJBVIs.

\begin{theorem}[Characterization of the value function] \label{uniq} Suppose that the above assumptions hold. The value function $u$ is then the unique viscosity solution of   the HJBVIs \eqref{edp} in the class of bounded continuous  functions,
 in the sense that $u$  is a viscosity sub and super solution of  \eqref{edp} with terminal condition $u(T,x)= g(x)$.

\end{theorem}
\noindent  Using the comparison principle (Lemma \ref{8.9}), we  show below some new results in the irregular case.
\section{Complementary results in the discontinuous case}\label{comp}

Using the previous results both in the discontinuous case and the continuous case, we show some additional properties of the value function when the terminal reward map is discontinuous.

Let $g$ be a Borelian function such that there exists $K>0$ and $p \in \mathbb{N}$ with $|g(x)| \leq K(1+|x|^p)$ and  satisfying $h_1(T,x) \leq g(x) \leq h_2(T,x)$ for all $x \in \mathbb{R}$. Let $g_*$ be the l.s.c. envelope of $g$.\\
We denote by $u^{g}(t,x)$ the value function of our problem associated with terminal reward $g$.

\begin{proposition}\label{suite}
Let $(g_n)_{n \in \mathbb{N}}$ be a non decreasing sequence of Lipschitz continuous maps satisfying for each  $n$, $h_1(T,x) \leq g_n(x) \leq h_2(T,x)$ and $|g_n(x)| \leq K(1+|x|^p)$, $x \in \mathbb{R}$, and such that $g_*= lim_{n \rightarrow \infty} \uparrow g_n$ (which exists by analysis results).  Then, we have
$$
u^{g_*}(t,x)= lim_{n \rightarrow \infty} \uparrow u^{g_n}(t,x).
$$
Moreover, for each $n$, $u^{g_n}$ is continuous and $u^{g_*}$ is l.s.c.
\end{proposition}
\dproof
By definition of the value function,
$$
u^{g_*}(t,x)= \sup_{\alpha \in \mathcal{A}_t^t} Y_{t,T}^{\alpha,t,x}[g_*(X_T^{\alpha,t,x})].
$$
By the continuity property of doubly reflected BSDEs with respect to the terminal condition (see Proposition \ref{continuity}) and by the comparison theorem, we have:
$$Y_{t,T}^{\alpha,t,x}[g_*(X_T^{\alpha,t,x})]=\lim_{n \rightarrow \infty} \uparrow Y_{t,T}^{\alpha,t,x}[g_n(X_T^{\alpha,t,x})].$$
From the two above equalities, we derive that
$$u^{g_*}(t,x)= \sup_{\alpha \in \mathcal{A}_t^t} \sup_{n \in \mathbb{N}} Y_{t,T}^{\alpha,t,x}[g_n(X_T^{\alpha,t,x})]=\sup_{n \in \mathbb{N}} \sup_{\alpha \in \mathcal{A}_t^t} Y_{t,T}^{\alpha,t,x}[g_n(X_T^{\alpha,t,x})]= \lim_{n \rightarrow \infty} \uparrow u^{g_n}(t,x).$$
By the continuity property of the value function in the continuous case (see Corollary \ref{certain}), $u^{g_n}$ is continuous because $g_n$ is continuous. Hence, $u^{g_*}$ is l.s.c. as supremum of continuous functions.
\fproof

\begin{corollary}\label{1}
We have:\\
1. If $g$ is l.s.c., then $u^{g}$ is l.s.c.\\
2. $u_{*}^{g} \geq u^{g_*}.$\\
3. The l.s.c. envelop $u_{*}^{g}$ of the value function $u^{g}$ is a viscosity supersolution of the HJBVIs \eqref{edp}, with terminal value greater than $g_*$.

\end{corollary}

\dproof
1. Since $g$ is l.s.c., $g_*=g$. Hence, Point 1. directly follows from the above Proposition.\\
2. As $g \geq g_*$, we have $u^{g} \geq u^{g_*}$. By Point 1, $u^{g_*}$ is l.s.c., and we thus obtain 
$u_{*}^{g} \geq u^{g_*}.$\\
3. By Point 2.\,, $u_{*}^{g}(T,x) \geq u^{g_*}(T,x)= {g_*}(x).$ Using Theorem \ref{existprob}, we derive that $u_{*}^{g}$ is 
a viscosity supersolution of the HJBVIs \eqref{edp}, with terminal value greater than $g_*$.
\fproof

\begin{corollary}
Suppose that ${\bf A}$ is compact. If $g$ is l.s.c., then $u^{g}$ is l.s.c.\,, and  is a viscosity supersolution of the HJBVIs \eqref{edp} with terminal condition $u^g(T,x)=g(x),\,\, x \in \mathbb{R}$.
\end{corollary}
\dproof 
By Point 1. of Corollary \ref{1}, as $g$ is l.s.c.,  $u^{g}$ is l.s.c. Hence, by Theorem \ref{existprob}, $u^g$
is a viscosity supersolution of  \eqref{edp}.
\fproof  

\begin{proposition}
Suppose that ${\bf A}$ is compact and  Assumption $\ref{9.8}$  holds. Let $v$ be a viscosity supersolution of the HJBVIs \eqref{edp} with $v(T,x) \geq g_*(x),\,\, x \in \mathbb{R}$. We then have $v \geq u^{g_*}.$ In other words, $u^{g_*}$ is the minimal viscosity supersolution of the HJBVIs \eqref{edp} with terminal value greater than $g_*$.
\end{proposition}

\dproof
Let $(g_n)_{n \in \mathbb{N}}$ be a non decreasing sequence of continuous maps satisfying 
the assumptions of Proposition \ref{suite}.
For each $n$, we have $v(T,x) \geq g_n(x),$ $x \in \mathbb{R}$. As $u^{g_n}$ is a continous viscosity solution, then it is a subsolution. We also have $u^{g_n}(T,x)=g_n(x)$, $x \in \mathbb{R}$. By using the comparison principle (see Lemma \ref{8.9}), we have $v \geq u^{g_n},$ for all $n$. Hence, we get
$$v \geq \lim_{n }\uparrow u^{g_n}=u^{g_*}, $$
where the last equality follows from Proposition \ref{suite}.
\fproof

\noindent By the above result, we have the following characterization of the value function when $g$ is l.s.c.

\begin{theorem}\label{5}
Suppose ${\bf A}$ is compact and  Assumption $\ref{9.8}$  holds. If $g$ is l.s.c., then the value function $u^g$ of our mixed problem  is l.s.c. and it is the minimal viscosity supersolution, with terminal value greater than $g$, of the HJBVIs \eqref{edp}.
\end{theorem}

\appendix
\section{Appendix}
\subsection{A decomposition property}
We show a property which ensures that under our assumptions (in particular $h_1$ or $h_2$  is ${\cal C}^{1,2}$),
then, for each $t \in [0,T]$ and for each $\alpha$ $\in$ ${\cal A}_t^t$, the barriers  $h_1(s,X_s^{\alpha,t,x}){\bf 1}_{s<T} + g(X^{\alpha,t,x}_T){\bf 1}_{s=T}$ and $h_2(s,X_s^{\alpha,t,x}){\bf 1}_{s<T} + g(X^{\alpha,t,x}_T){\bf 1}_{s=T}$ 
 satisfy Mokobodzki's condition (see \eqref{mo}).

\begin{proposition}\label{Moko}
Let $h :[0,T] \times \mathbb{R} \rightarrow \mathbb{R}$  be ${\cal C}^{1,2}$ with bounded derivatives,
and let $g : \mathbb{R} \rightarrow \mathbb{R}$ be a Borelian map. 
Suppose that 
$  |h(t,x)| +|g(x)| \leq C(1+ |x|^p).$ 
Let $t \in [0,T]$ and  $\alpha$ $\in$ ${\cal A}_t^t$. 
 There exist  two  nonnegative  ${\mathbb F}^t$-supermartingales $H^{\alpha,t,x}$ and $H^{' \alpha,t,x}$ in ${\cal S}_t^2$ such that 
 \begin{equation*}
  h(s,X_s^{\alpha,t,x}){\bf 1}_{s<T} + g(X^{\alpha,t,x}_T){\bf 1}_{s=T}= H^{\alpha,t,x}_s - H^{' \alpha,t,x}_s , 
 \,\,\,t \leq s \leq T \quad {\rm a.s.}
 \end{equation*}
\end{proposition}

 \dproof 
 Set  $\xi_s :=h(s,X_s^{\alpha,t,x}){\bf 1}_{s<T} + g(X^{\alpha,t,x}_T){\bf 1}_{s=T}$.
 By applying It\^o's formula to $h(s,X_s^{\alpha,t,x})$, one can show that there exist three ${\mathbb F}^t$-predictable processes $(f_r)$,
 $(\varphi^1_r)$ and $ (\varphi^2 (r,.))$  such that
 $dh(r,X_r^{\alpha,t,x})= f_rdr +\varphi^1_r dW^t_r +\int_E  \varphi^2 (r,e) d\tilde N^t(dr,de)$, 
 with $|f_r| + | \varphi^1_r| + || \varphi^2 (r,.)||_{\nu} 
 \leq K(1 + |X_{r^-}^{\alpha,t,x}|^p)$. \\
  Integrating this 
 equation between $s \geq t$ and $T$, and then taking the conditional expectation with respect to 
 ${\cal F}^t_s$, we derive that
 $$h(s,X_s^{\alpha,t,x})= E[h(T,X_T^{\alpha,t,x}) + \int_s ^T f_r dr \,| \, {\cal F}^t_s]= I_s - I'_s, $$
 where $(I_s)$ and $(I'_s)$ are defined for each $s \in [t, T]$ by
 $I_s := E[h^+_1(T,X_T^{\alpha,t,x}) + \int_s ^T f^+_r dr  \,| \, {\cal F}^t_s]$ and 
 $I'_s := E[h^-_1(T,X_T^{\alpha,t,x}) + \int_s ^T f^-_r dr  \,| \, {\cal F}^t_s]$. They are both nonnegative
  supermartingales belonging to ${\cal S}_t^2$. For each $s \in [t,T]$, we have
  $h(s,X_s^{\alpha,t,x})=I_s - I'_s$ . 
  For each $s \in [t,T]$, set
\begin{equation}\label{ct}
\Tilde{\xi}_s:=\xi_s-{E}[g(X^{\alpha,t,x}_T)| {\cal F}^t_s]. 
\end{equation}
Let us now show that  there exist two nonnegative supermartingales $\tilde H$ and $\tilde H'$ such that
 \begin{equation}\label{tildem}
\Tilde{H}_T=\Tilde{H'}_T=0\quad {\rm and} \quad \Tilde{\xi}_s   = \tilde H_s - \tilde H'_s, \,\, t \leq s \leq T.\end{equation}
By \eqref{ct}, we have
$\Tilde{\xi}_T = 0 $. 
We thus have that for each $s \in [t,T]$,
\begin{equation*}
\Tilde{\xi}_s=(h(s,X_s^{\alpha,t,x})-{E}[g(X^{\alpha,t,x}_T)| {\cal F}^t_s] ) {\bf 1}_{s<T}
= (I_s - I'_s-{E}[g(X^{\alpha,t,x}_T)| {\cal F}^t_s] ) {\bf 1}_{s<T}.
\end{equation*}
Now, the processes $(I_s + {E}[g^-(X^{\alpha,t,x}_T)| {\cal F}^t_s])_{t \leq s \leq T}$ and 
$(I'_s + {E}[g^+(X^{\alpha,t,x}_T)| {\cal F}^t_s])_{t \leq s \leq T}$ are nonnegative supermartingales  as the sum of two  nonnegative supermartingales.\\
It follows that  the processes $(\tilde H_s)_{t \leq s \leq T}$ and 
$(\tilde H'_s)_{t \leq s \leq T}$ defined by\\
 $\tilde H_s:= (I_s + {E}[g^-(X^{\alpha,t,x}_T)| {\cal F}^t_s]){\bf 1}_{s<T}$ and
$\tilde H'_s:= (I'_s + {E}[g^+(X^{\alpha,t,x}_T)| {\cal F}^t_s]){\bf 1}_{s<T}$, $t \leq s \leq T$, 
are also nonnegative supermartingales, which are moreover equal to $0$ at time $T$. They satisfy the equality 
$\Tilde{\xi}_s= \tilde H_s - \tilde H'_s$ for each $s \in [t,T]$. 
Hence, $\tilde H$ and $\tilde H'$ satisfy equalities \eqref{tildem}. \\
 Using equality \eqref{ct}, we derive that the processes $H$ and $H'$ defined by 
 $H_s:= \tilde H_s + {E}[g^+(X^{\alpha,t,x}_T)| {\cal F}^t_s]$ and $H'_s:= \tilde H'_s + {E}[g^-(X^{\alpha,t,x}_T)| {\cal F}^t_s]$ are nonnegative supermartingales satisfying
 $ {\xi}_s =  H_s - H'_s,$ $t \leq s \leq T$, which ends the proof. 
 \fproof
 
\noindent Note that we cannot apply
 It\^o's formula to ${\bar h}(s,X_s^{\alpha,t,x})$ with 
 ${\bar h}( s,x):= h( s,x){\bf 1}_{s<T}  + g(x){\bf 1}_{s=T}$ since 
 $g$ is irregular.
The change of variable \eqref{ct}
 allows us to deal more easily with the irregularity of $g$.

\subsection{Measurability and splitting properties}\label{A1}
For each $\omega \in \Omega$, let $^s\omega:= (\omega_{r \wedge s})_{0 \leq r \leq T}$. 
Let $S^s$ (resp. $T^s$) be the operator defined on $\Omega$ by 
$S^s( \omega) := $$^s\omega$ (resp. $T^s( \omega) := \omega^s$). 
Note that $S^s$ and $T^s$ are independent and for each $ \omega \in \Omega$  we have
$\omega = S^s(\omega) + T^s(\omega) {\bf 1}_{]s,T]}=$ $ ^s\omega +  \omega^s {\bf 1}_{]s,T]}$. 
When there is no ambiguity, we identify $\omega$ with $(^s\omega, \omega^s).$

Let $t \in [0,T]$, $\alpha$ $\in$ ${\cal A}_t^t$ and
$s \geq t$. For $P$-almost every $ \omega$, the process $\alpha(^s \omega, T^s)$ (denoted also by $\alpha (^s \omega, .)$) 
defined by
$$\alpha(^s \omega, T^s): \Omega  \times [s,T] \rightarrow {\mathbb R} \,; \,(\omega ', r) \mapsto \alpha_r(^s \omega, T^s(\omega '))$$
belongs to  
$ \mathcal{A}_{s}^{s}$ (see Lemma 3.3 in \cite{DQS3}).

Recall that the forward process satisfies the splitting property (see \cite{DQS3}):
for all $t \in [0,T]$, 
$\alpha$ $\in$ ${\cal A}_t^t$ and
$s \geq t$, for almost every $\omega \in \Omega$, setting $\tilde \omega= $ $^{s}\omega$, the process 
 $X^{\alpha,t,x} (\tilde \omega, T^s)$ (denoted also by $X^{\alpha,t,x} (\tilde \omega, .)$) coincides with the solution of the forward SDE 
 on $\Omega \times[s,T]$, driven by $ W^s$ and $\Tilde{N}^s$, associated with control 
 $ (\alpha_r(\tilde \omega,\cdot))_{r \geq s}$ and  initial condition $X_s^{\alpha(\tilde \omega),t,x}(\tilde \omega)$ at time $s$.
 
 By similar arguments as in \cite{DQS3}, we have an analogous property for doubly reflected BSDEs: 

\begin{proposition}\label{egal}(Splitting property for doubly reflected BSDEs)

Let $t \in [0,T]$, $\alpha$ $\in$ ${\cal A}_t^t$ and
$s \in [t,T]$. 
%
 For almost every $\omega \in \Omega$, setting $\tilde \omega= $ $^{s}\omega$, the process 
 $(Y_r^{\alpha,t,x} (\tilde \omega, .)) _{ s \leq r \leq T}$ coincides with the solution of the doubly reflected BSDE 
 on $\Omega \times[s,T]$, associated with driver 
$f^{\alpha(\tilde \omega,.),s,X_s^{\alpha(\tilde \omega),t,x}(\tilde \omega)}$, barriers
 $ h_i(r, X_r^{\alpha(\tilde \omega,T^s),s,X_s^{\alpha(\tilde \omega),t,x}(\tilde \omega)})_{r < T}$, $i=1,2$, terminal condition $ g(X_T^{\alpha(\tilde \omega,T^s),s,X_s^{\alpha(\tilde \omega),t,x}(\tilde \omega)})$, filtration $\mathbb{F}^s$,  and driven by $ W^s$ and $\Tilde{N}^s$. In particular, we have
 \begin{equation}\label{abcd}
 Y_s^{\alpha,t,x} (\tilde \omega)  = Y_{s}^{ \alpha(\tilde\omega, \cdot),s,X_s^{\alpha(\tilde \omega),t,x}
 (\tilde \omega)} = u^{ \alpha(\tilde\omega, \cdot)}(s,X_s^{\alpha(\tilde \omega),t,x}(\tilde \omega)).
\end{equation}

\end{proposition}


Let $\eta \in {L}^2(\mathcal{F}_{s}^t).$ 
Since $\eta$ is $\mathcal{F}_{s}$-measurable, up to a $P$-null set, it can be written as a 
measurable map, still denoted by $\eta$, of the past trajectory $^s\omega$.
For each $\omega \in {\Omega}$, by using the definition of the function $u$, we have:
\begin{equation}\label{aux}
u(s,\eta(^s \omega))=\sup_{\alpha \in \mathcal{A}_s^s} u^{\alpha}(s, \eta(^s\omega)).
\end{equation}

 For each $(t,s)$ with $s \geq t$, we introduce the set $\mathcal{A}_s^t$ of restrictions to $[s,T]$ of the controls in $\mathcal{A}_t^t$. They can also be identified to the controls $\alpha$ in $\mathcal{A}_t^t$ which are equal to $0$ on $[t,s]$. \\
 Using the measurability and continuity results of Section 2.1, we show the following result. 
 
 \begin{proposition} \label{mesu2} 
Let $t \in [0,T]$, $s \in [t,T]$ and $\eta \in {L}^2(\mathcal{F}_s^t).$ Let $\varepsilon >0$. 
 Suppose that ${\bf A}$ is compact. There exists 
$\alpha^{\varepsilon}$ $\in \mathcal{A}_s^t$ such that,
for almost every $ \omega \in \,  \Omega$,  
 \begin{equation}\label{optimum}
u(s,\eta(^s \omega)) \,\leq \, u^{\alpha^{\varepsilon}(^s\omega,.)}(s, \eta(^s\omega))+\varepsilon.
\end{equation}
Moreover, the $\varepsilon$-optimal control $\alpha^{\varepsilon}$  can be constructed so that it
depends on the past trajectory $^s \omega$
  only through $\eta(^s \omega)$.\\
If 
Assumption \ref{H2} holds, this result still holds when
${\bf A}$  is a nonempty closed subset of ${\mathbb R}$, not necessarily compact. 
\end{proposition} 

\dproof 
To simplify notation, we suppose that $t=0$. We introduce  $\Omega^s$  the set of the restrictions to $[s,T]$ of the paths $\omega$ $\in$ $\Omega$. 
By Lemma \ref{mesu},  the map $(x, \alpha) \mapsto u^{ \alpha}(s,x)$ is  
$ \mathcal{B}(\mathbb{R})\otimes \mathcal{B}(\mathcal{A}_s^s )$-measurable. 

By Proposition 7.50 in \cite{BS} together with  a result of 
measure theory (see e.g. Lemma A.3 in \cite{DQS3}),
we derive that there exists a Borelian map  ${\hat \alpha}^{\varepsilon}:$ 
$\mathbb{R} \rightarrow \mathcal{A}_s^s$\,; $\,x\mapsto {\hat \alpha}^{\varepsilon}(x, \cdot)$ 
 such that  
\begin{equation}\label{n}
u(s,x) \leq u^{ {\hat \alpha}^{\varepsilon}(x,\cdot)}(s, x)+\varepsilon  \quad {\rm for} \quad {\rm Q}-{\rm almost\,\, every}\,\, x \in \mathbb{R},
\end{equation} 
where $Q$ is the law of $\eta$ under $P$.
\noindent 
Let $\{e^{i}, i \in \mathbb{N} \}$ be a countable orthonormal basis of the separable Hilbert space $\mathbb{H}_s^2$. 
For each $x$, the map ${\hat \alpha}^{\varepsilon}_u(x, \cdot)$, which belongs to $\mathbb{H}_s^2$, admits the following decomposition: 
${\hat \alpha}^{\varepsilon}_u(x, \cdot)=\sum_{i}\beta^{i, \varepsilon}(x)e^i_u(\cdot)$ $dP\otimes du$-a.s. and in $\mathbb{H}_s^2$, where $\beta^{i,\varepsilon}: \mathbb{R} \mapsto \mathbb{R}$ is a Borelian map defined  by $\beta^{i,\varepsilon}(x)=<{\hat \alpha}^{\varepsilon}(x,\cdot),e^{i}(\cdot)>_{\mathbb{H}_s^2}$. 
Consider the predictable process $\alpha^{\varepsilon}$  defined for each  $(r,\omega) \in [s,T] \times \Omega$  by
$\alpha^{\varepsilon}_r(\omega):=   \sum_{i}\beta^{i, \varepsilon}(\eta(^s\omega)) e^i_r( \omega^s)$.
Since ${\bf A}$ is bounded, $\alpha^{\varepsilon}$ is bounded, and hence  belongs to $\mathbb{H}_s^2$ and thus to ${\cal A}_s$, which gives the desired result.

When ${\bf A}$ is only a nonempty closed subset of ${\mathbb R}$, not necessarily bounded,  the control
 $\hat \alpha^{\varepsilon}$ can also be supposed to be bounded because 
 the set of bounded controls in ${\cal A}_s ^s$, denoted by  ${\cal A}_s ^{s,b}$,  is dense in ${\cal A}_s ^s$. 
Indeed, let $a \in {\bf A}$. If $\alpha$ $\in {\cal A}_s ^s$, the  bounded controls $\alpha_r {\bf 1}_{| \alpha_r | \leq n   }+ a {\bf 1}_{| \alpha_r| > n   }$ belong to $ \mathcal{A}_s ^s$ and converge to $\alpha$ in $\mathbb{H}_s^2$ as $n$ tends to $+ \infty$.\\
 Now, by the continuity Assumption \ref{H2}, 
$u^{\alpha}$ is continuous 
with respect to $\alpha$ (see Lemma 
\ref{mesura}). Hence, 
 $$u(s,x) =\sup_{\alpha \in \mathcal{A}^s_s} u^{\alpha}(s,x)= \sup_{\alpha \in {\cal A}_s ^{s,b} } u^{\alpha}(s,x).$$
 It follows that there exists a bounded Borelian map  ${\hat \alpha}^{\varepsilon}:$ 
$\mathbb{R} \rightarrow \mathcal{A}_s^{s}$\,; $\,x\mapsto {\hat \alpha}^{\varepsilon}(x, \cdot)$ 
which satisfies  inequality \eqref{n}.
The proof is thus complete.
\fproof


\noindent By applying this property to $\varepsilon = \frac{1}{n}$ for each $n\in\mathbb{N}$, 
we derive that under Assumption \ref{H2}, there exists an ``optimizing sequence"  for the value function $u(s, \eta(\cdot))$.

\begin{corollary}\label{selectbis} Suppose 
  Assumption \ref{H2} holds. 
Let $t \in [0,T]$, $s \in [t,T]$ and $\eta \in {L}^2(\mathcal{F}_s^t).$  
 There exists a sequence $(\alpha^n)_{n\in\mathbb{N}} \in \mathcal{A}^t_s$ such that 
 for almost every $\omega \in \Omega$,
 $$u(s, \eta(^s \omega))=\lim_{n \rightarrow \infty}u^{\alpha^n(^s \omega,\cdot)}
 (s,\eta(^s \omega)).$$
 Moreover, the processes $\alpha^n$, $n\in\mathbb{N}$, can be chosen so that they
depend on the past trajectory $^s \omega$
  only through $\eta(^s \omega)$.\\
\end{corollary}


\subsection{Complementary results on BSDEs and reflected BSDEs with jumps}
We first state a version of Fatou lemma for ${\cal E}^f$-conditional expectations  (or equivalently for BSDEs) where the limit involves both terminal condition and terminal time.
\begin{lemma}[A Fatou lemma for BSDEs]\label{fatou1}
Let $T>0$.  Let $f$ be a given Lipschitz driver  satisfying Assumption  \ref{Royer}. 
Let $(\theta^n)_{ n \in \mathbb{N}} $ be a non increasing sequence of stopping times in 
$\mathcal{T}$,  converging a.s. to  $\theta \in \mathcal{T}$ as 
$n$ tends to $\infty$. Let   $(\xi^{n})_{ n \in \mathbb{N}} $ be a sequence of random variables such that 
$ \mathbb{E}[\sup_{n}( \xi^n)^2]< + \infty$, and  for each $n$,
$\xi^{n}$ is ${\mathcal F}_{\theta^{n}}$-measurable.
Then,  for each  stopping time $\tau $ with $\tau \leq \theta$ a.s.\,,  we have
 $${\cal E}^f_{\tau,\theta}(\liminf_{n \rightarrow + \infty} \xi^{n} ) \leq \liminf_{n \rightarrow + \infty} {\cal E}^f_{\tau,\theta^{n}}(\xi^n )  \quad {\rm and} \quad {\cal E}^f_{\tau,\theta}(\limsup_{n \rightarrow + \infty} \xi^{n} ) \geq \limsup_{n \rightarrow + \infty} {\cal E}^f_{\tau,\theta^{n}}(\xi^n )  \quad a.s.\,$$
  \end{lemma}

\dproof 
Let us show the first inequality.
For each $n$, by the nondecreasing property of ${\cal E}^f$, we have 
${\cal E}^f_{\tau,\theta^{n}} ( \inf_{p \geq n}  \xi^{p}) \leq {\cal E}^f_{\tau,\theta^{n}} (\xi^{n})$ a.s.\, 
Thus, 
 $ \liminf_{n \rightarrow + \infty}{\cal E}^f_{\tau,\theta^{n}} ( \inf_{p \geq n}   \xi^{p}) \leq \liminf_{n \rightarrow + \infty}{\cal E}^f_{\tau,\theta^{n}} (\xi^{n})$ a.s.\,
Since 
 $\lim_{n \rightarrow + \infty} \inf_{p \geq n}   \xi^{p} = \liminf_{n \rightarrow + \infty} \xi^{n}$ a.s.\,, by a continuity property the ${\cal E}^f$-conditional expectation (see Proposition A.6 in \cite{16}), we get 
$ {\cal E}^f_{\tau,\theta}(\liminf_{n \rightarrow + \infty} \xi^{n} ) = \lim_{n \rightarrow + \infty}{\cal E}^f_{\tau,\theta^{n}} ( \inf_{p \geq n} \xi^{p})$ a.s.\,, from which the desired inequality follows. The proof of the second one is similar.
\fproof

\noindent Using this lemma, we now show a continuity property of the solutions of reflected BSDEs with respect to both terminal time and terminal condition, which extends the result established in \cite{DQS3} at time $0$ (see Lemma 3.10) to any  time 
$\tau \in \mathcal{T}$.
\begin{lemma}[A continuity property for reflected BSDEs] \label{CRBSDE}
Let $T>0$.  Let $(  \eta_t)$ be an RCLL process in $\mathcal{S}^2.$ Let $f$ be a given Lipschitz driver satisfying Assumption \ref{Royer}. 
Let $(\theta^n)_{ n \in \mathbb{N}} $ be a non increasing sequence of 
stopping times in 
$\mathcal{T}$,  converging a.s. to  $\theta \in \mathcal{T}$ as 
$n$ tends to $\infty$. Let   $(\xi^{n})_{ n \in \mathbb{N}} $ be a sequence of random variables such that 
$ \mathbb{E}[{\rm ess }\sup_{n}( \xi^n)^2]< + \infty$, and  for each $n$,
$\xi^{n}$ is ${\mathcal F}_{\theta^{n}}$-measurable. Suppose  that $\xi^{n}$ converges a.s. to an ${\mathcal F}_{\theta}$-measurable random variable $\xi$ as 
$n$ tends to $\infty$. Suppose that 
\begin{equation}\label{assi}
  \eta_{\theta} \leq \xi \quad {\rm a.s.}
\end{equation}
Let $Y_{.,\theta^{n}}(\xi^n )$; $Y_{.,\theta}(\xi )$ be the solutions of the reflected BSDEs associated with driver $f$, obstacle $(  \eta_s)_{ s < \theta^n}$ (resp. $(  \eta_s)_{ s < \theta}$), terminal time $\theta^n$ (resp. $\theta$), terminal condition $\xi^n$ (resp. $\xi$).
Then,  for each  stopping time $\tau$ with $\tau \leq \theta$ a.s.\,,  
 $$Y_{\tau,\theta}(\xi ) = \lim_{n \rightarrow + \infty} Y_{\tau,\theta^{n}}(\xi^n ) \quad a.s.$$
When for each $n$,  $\theta_n= \theta$ a.s.\,, the result still holds without Assumption \eqref{assi}.

 \end{lemma}

 \dproof  Note that the case where $\tau=0$ has been solved in \cite{DQS3} by using  the classical a priori estimates on reflected BSDEs. The case where $\tau$ is 
 any stopping time requires some additional arguments. It could be shown by using again It\^o's calculus. We adopt here a less classical approach which requires less computations.\\
  {\bf Step 1}: Let us first consider the simpler case when for each $n$, $\theta_n = \theta$ a.s.\\
  Using the a priori estimates on reflected BSDEs provided \cite{DQS} (see Proposition A.1 in \cite{DQS}) and
  the convergence of $\xi^ n$ to $\xi$, one can show that  $Y_{\tau,\theta}(\xi ) = \lim_{n \rightarrow + \infty} Y_{\tau,\theta}(\xi^n ) $ a.s.\\ 
{\bf Step 2}:  Let us now consider the general case. The difficulty is here to deal with the variation of 
the terminal time together with the presence of the obstacle $(  \eta_t)$. In particular, Proposition A.1 in \cite{DQS} is not appropriate to this case.
 By the flow property for reflected BSDEs, we have:
  \begin{align*}
 Y_{\tau, \theta_n}(\xi^n)=Y_{\tau, \theta}(Y_{\theta, \theta_n}(\xi^n)).
 \end{align*}
 By step 1,  it is thus sufficient to show that $\lim_{n \rightarrow} Y_{\theta, \theta_n}(\xi^n)=\xi$ a.s.\\
 Since the solution of the reflected BSDE associated with terminal condition  $\xi^n$ and terminal time $\theta_n$ is greater than the solution of the nonreflected BSDE associated with terminal condition  $\xi^n$ and terminal time $\theta_n$, we have:
 \begin{equation}\label{ineq-}
 Y_{\theta, \theta_n}(\xi^n) \geq \mathcal{E}_{\theta, \theta_n}(\xi^n)  \text{ a.s. }
 \end{equation}
  Now, by the continuity property of BSDEs with respect to both terminal time and terminal condition (Proposition A.6 in \cite{16}), we have $\lim_{n \rightarrow \infty} \mathcal{E}_{\theta, \theta_n}(\xi^n)= \mathcal{E}_{\theta, \theta}(\xi)=
 \xi$ a.s. Taking the $\lim \inf$ in $\eqref{ineq-}$, we obtain:
 \begin{align*}
 \lim \inf_{n \rightarrow \infty} Y_{\theta, \theta_n}(\xi^n) \geq \lim_{n \rightarrow \infty} \mathcal{E}_{\theta, \theta_n}(\xi^n)= \xi, \,\,\, \text{ a.s.}
 \end{align*}
  It remains to show that:
 \begin{equation}\label{resultat}
  \lim \sup_{n \rightarrow \infty} Y_{\theta, \theta_n}(\xi^n) \leq \xi, \,\,\, \text{ a.s.}
 \end{equation}
 By the characterization of the solution of a reflected BSDE, we obtain:
 \begin{align*}
 Y_{\theta, \theta^n}= ess\sup_{\tau \in \mathcal{T}_{\theta}} \mathcal{E}_{\theta, \tau \wedge \theta_n}(  \eta_\tau \textbf{1}_{\tau < \theta_n}+ \xi^n \textbf{1}_{\tau \geq \theta_n}).
 \end{align*}
 Fix $\varepsilon>0$. By the second assertion of Theorem 3.3 in \cite{17}, there  exists an $\varepsilon$-optimal stopping time $\tau_n^\varepsilon \in \mathcal{T}_\theta$, 
 that is  such that
 \begin{equation}\label{ineq1-}
 Y_{\theta, \theta_n}(\xi^n) \leq \mathcal{E}_{\theta, \tau_n^\varepsilon \wedge \theta_n}(  \eta_{\tau_n^\varepsilon} \textbf{1}_{\tau_n^\varepsilon < \theta_n}+ \xi^n \textbf{1}_{\tau_n^\varepsilon \geq \theta_n}) + \varepsilon \quad {\rm a.s.}
 \end{equation}
 Note that we have the following property: let $X,X', X_n, X'_n$, $n \in \mathbb{N}$,
   be real valued random variables 
 with $X \leq X'$, and let $A_n$, $n \in \mathbb{N}$ be measurable sets of $\mathcal{F}_T.$
  \begin{equation}\label{ineq6}
  {\rm If } \, X_n \rightarrow X  \, {\rm and } \, X'_n \rightarrow X'  \, \text{  a.s., then } \, \, \lim \sup_{n \rightarrow \infty} (X_n \textbf{1}_{A_n}+ X'_n \textbf{1}_{A_n^c}) \leq X' \,  {\rm a.s.}
  \end{equation}
Now, for each $n$, $\tau_n^\varepsilon \wedge \theta_n \geq \theta$ a.s. and $\tau_n^\varepsilon \wedge \theta_n $ tends to $\theta$ a.s. as $n \rightarrow + \infty$.
Hence, by the right-continuity property of the obstacle $(  \eta_t)$, we get $  \eta_{\tau_n^\varepsilon \wedge \theta_n} \rightarrow   \eta_{\theta} \leq \xi$ a.s.\,, where the last inequality holds by Assumption \ref{assi}.
 By applying Property \eqref{ineq6} and since $\xi^n \rightarrow \xi$ a.s., we thus obtain
 \begin{align*}
 \lim \sup_{n \rightarrow \infty} (  \eta_{\tau_n^\varepsilon} \textbf{1}_{\tau_n^\varepsilon < \theta_n}+\xi^n \textbf{1}_{\tau_n^\varepsilon \geq \theta_n}) \leq \xi \text{ a.s.}
 \end{align*}
 Now, the Fatou property for BSDEs (see Lemma \ref{fatou1}) together with \eqref{ineq1-} implies that\\ $\lim \sup_{n \rightarrow \infty} Y_{\theta, \theta_n}(\xi^n) \leq \xi+\varepsilon$ a.s.
 Since this inequality holds for each $\varepsilon>0$, we derive inequality \eqref{resultat}. The proof is thus complete.
  \fproof


\small
\begin{thebibliography}{15}

 \bibitem{ALM} Alario-Nazaret, M. Lepeltier, J.P. and Marchal, B. (1982). Dynkin games. (Bad Honnef Workshop on stochastic processes), {\em Lecture Notes in control and Information Sciences} 43, Springer-Verlag, Berlin.



%










\bibitem{BH} Bayraktar E. and Huang, Y.J., On the multi-dimensional controller-and-stopper game,  {\em 
SIAM Journal of Control and Optimization} (2013), 51, 1263-1297.



\bibitem{BL} Bensoussan, A. and Lions J-L, {\em Applications des in\'equations variationnelles en contr\^ole stochastique} Volume 6 de M\'ethodes math\'ematiques de l'information, Dunod, 1988.


\bibitem{BF} Bensoussan, A. and Friedman, A. 
Non-linear variational inequalities and differential games with stopping times, {\em J. Funct. Anal.}(1974), 16, pp. 305-352.

\bibitem{BS} Bertsekas D. and  Shreve S., 
{\em Stochastic Optimal Control: The Discrete Time Case}, Academic Press, Orlando (1978).

\bibitem{Bismut}   Bismut J.M., Sur un probl\`eme de Dynkin, {\em Z.Wahrsch. Verw. Gebiete} 39 (1977) 31--53.


\bibitem{BT} Bouchard, B. and N. Touzi. Weak Dynamic Programming Principle for Viscosity Solutions,  {\em SIAM J Control Optim}, 2011, 49 (3), 948-962.


%



\bibitem{Buckdahn1} Buckdahn, R. and Nie, T. (2014), Generalized Hamilton-Jacobi-Bellman equations with Dirichlet boundary and stochastic exit time optimal control problem, \\
http://arxiv.org/pdf/1412.0730v4.pdf
\bibitem{BuckdahnR} Buckdahn, R. and Li, J. (2010), Stochastic Differential Games with Reflection and Related Obstacle Problems for Isaacs Equations, Acta Mathematicae Applicatae Sinica, English Series.
\bibitem{Buckdahn2} Buckdahn, R. and Li, J. (2009), Probabilistic interpretation for systems of Isaacs equations with two reflecting barriers, {\em Nonlinear Differ.Equ. Appl. }16, pp. 381-420.






\bibitem{CCP}  Choukroun S., Cosso A. and Pham H.,
Relected BSDEs with nonpositive jumps, and controller-and-stopper games,
 {\em Stochastics Processes and their Applications}, vol. 125(2), 597-633, 2015.


\bibitem{CM} Cr\'epey, S. and Matoussi, A., Reflected and Doubly Reflected BSDEs with jumps, {\em Annals of App. Prob.} 18(5),  2041-2069 (2008).

\bibitem{CK} Cvitani\'c J. and Karatzas, I.\,, Backward stochastic differential equations with reflection and Dynkin games, \textit{Annals of Prob.} {\bf 1996}. \textbf{24}, n.4 2024-2056.











\bibitem{DM1}
Dellacherie, C. and  Meyer, P.-A. (1975).
 \textit{Probabilit\'es et Potentiel, Chap. I-IV}. Nouvelle \'edition. Hermann. {\bf MR}{0488194}
\bibitem{DQS}  Dumitrescu, R., Quenez M.C., Sulem A., Optimal stopping for dynamic risk measures with jumps and obstacle problems, {\em Journal of Optimization Theory and Applications},  2014, DOI:10.1007/s10957-014-0635-2.

\bibitem{DQS2}  Dumitrescu, R., Quenez M.C., Sulem A., (2013) Generalized Dynkin Games and Doubly reflected BSDEs with jumps, to appear in {\em Electronical Journal of Probability}, http://arxiv.org/pdf/1310.2764v2.pdf

\bibitem{DQS3}  Dumitrescu, R., Quenez M.C., Sulem A.,  A Weak Dynamic Programming Principle for Combined Optimal Stopping/ Stochastic Control with ${\cal E}^f$-expectations, to appear in {\em SIAM Journal on Control and Optimization}, http://arxiv.org/pdf/1407.0416v4.pdf 

\bibitem{DQS4}  Dumitrescu, R., Quenez M.C., Sulem A., Game options in an imperfect market with default, arXiv:1511.09041, November 2015






%

%
%


%







%
%


\bibitem{H} Hamad\`ene, S., Mixed zero-sum stochastic differential game and American game options, {\em SIAM J. Control Optim.}, 45(2), 
(2006), 496-518.


\bibitem{HH} Hamad\`ene, S. and  Hassani, M.,  BSDEs with two reacting  barriers driven by a Brownian motion
and an independent Poisson noise and related Dynkin game, {\em Electron. J. Probab.},
 11(5) (2006),  121-145.


%
\bibitem{HL} Hamad\`ene, S. and Lepeltier, J.-P.\,, Reflected BSDEs and mixed game problem, {\em Stochastic Process. Appl.} 85(2000) 177-188.

\bibitem{HO} Hamad\`ene, S.  and Ouknine, Y., Reflecting  Backward SDEs  with general jumps, {\em Teor. Veroyatnost. i Primenen.}  (2015), 60(2),  357-376.

\bibitem{HW} Hamad\`ene, S. and Wang, H. : BSDEs with two RCLL reflecting obstacles driven by a Brownian motion 
and  Poisson random measure and related mixed zero-sum games, {\em Stochastic Process. Appl.} 119 (2009) 2881--2912.

\bibitem{J}   Jacod, J. (1979).
 \textit{Calcul Stochastique et Probl\`emes de martingales}, Springer.

%
\bibitem{KZ} Karatzas I. and  Zamfirescu I., Martingale approach to stochastic differential games of control and stopping,  {\em Annals of Probability}, 36 (4) 1495-1527, 2008.

\bibitem{Kifer} Kifer. Y. and  Yu I.,  Game options,  {\em Finance and Stochastics}, (4) 443-463, 2012.


\bibitem{KiferYu} Kifer Y., Dynkin Games and Israeli options,   {\em ISRN Probability and Statistics}, ID856458, 2013.


\bibitem{KQC} Kobylanski M., Quenez M.-C. and Roger de Campagnolle M., Dynkin games in a general framework, {\em Stochastics}, 2013.



\bibitem{OS} {\O }ksendal, B. and Sulem, A. {\em Applied Stochastic Control of Jump Diffusions}, 2nd Ed., Universitext, Springer, Berlin, 2007.

 

\bibitem{Peng92}  Peng, S.,
A generalized dynamic programming principle and Hamilton-Jacobi-Bellman-Equation, 
 {\em Stochastics and Stochastics Reports}, 38, 1992.









\bibitem{16} Quenez M.-C. and Sulem A., BSDEs with jumps, optimization and applications to dynamic risk measures, {\em Stochastic Processes and Applications} 123 (2013) 3328-3357.

\bibitem{17} Quenez M.-C.  and  Sulem A., Reflected BSDEs and robust optimal stopping for dynamic risk measures with jumps, {\em Stochastic Processes and Applications} 124, (2014) 3031--3054 .



%
%
%
%









\end{thebibliography}
\end{document}